\documentclass{smjour}

\usepackage{combdesi}
\usepackage{amsmath}
\usepackage{graphicx}
\usepackage{longtable}
\usepackage{url}
\usepackage{enumerate}

\newcommand{\A}{{\mathcal A}}
\newcommand{\B}{{\mathcal B}}

\newcommand{\OO}{{\mathcal O}}
\newcommand{\SSS}{{\mathcal S}}

\newcommand{\vu}{{\sf u}}

\newcommand{\2}{{\underline{2}}}
\newcommand{\3}{{\underline{3}}}
\newcommand{\4}{{\underline{4}}}
\newcommand{\5}{{\underline{5}}}
\newcommand{\6}{{\underline{6}}}

\begin{document}

\yearofpublication{}
\volume{}
\cccline{}

\received{}
\accepted{}

\authorrunninghead{CHEE and KASKI}
\titlerunninghead{An Enumeration of Graphical Designs}

\title{An Enumeration of Graphical Designs}

\author{Yeow Meng Chee$^{1,2}$, Petteri Kaski$^3$}

\affil{\\ \footnotemark[1] Division of Mathematical Sciences,
School of Physical and Mathematical Sciences,
Nanyang Technological University,
Singapore 637616.}

\affil{\\ \footnotemark[2] Card View Pte. Ltd.,
41 Science Park Road,
\#04-08A The Gemini,
Singapore Science Park II,
Singapore 117610.}

\affil{\\ \footnotemark[3] Helsinki Institute for Information Technology HIIT,
Department of Computer Science,
University of Helsinki,
P.O. Box 68, 00014 University of Helsinki,
Finland.}

\abstract{
\boldmath
Let $\Psi(t,k)$ denote the set of pairs $(v,\lambda)$ for which there exists a
graphical $t$-$(v,k,\lambda)$ design. Most results on graphical designs have
gone to show the finiteness of $\Psi(t,k)$ when $t$ and $k$ satisfy certain conditions.
The exact determination of $\Psi(t,k)$ for specified $t$ and $k$ is a hard problem and
only $\Psi(2,3)$, $\Psi(2,4)$, $\Psi(3,4)$, $\Psi(4,5)$, and $\Psi(5,6)$ have been determined.
In this paper, we determine completely the sets $\Psi(2,5)$ and $\Psi(3,5)$.
As a result, we find more than 270000 inequivalent graphical designs, and
more than 8000 new parameter sets for which there exists a graphical design.
Prior to this,
graphical designs are known for only 574 parameter sets.
}

\begin{article}

\section{Introduction}
For a finite set $X$ and a nonnegative integer $t$, the set of all $t$-subsets of $X$
is denoted $\binom{X}{t}$.
A {\em $k$-uniform set system} is a pair $(X,\B)$, where $X$ is a finite set of elements called {\em points}
and $\B\subseteq\binom{X}{k}$. Elements of $\B$ are called {\em blocks}. The {\em order} of
$(X,\B)$ is the number of points, $|X|$.
A {\em design} with parameters $t$-$(v,k,\lambda)$ is a
$k$-uniform set system $(X,\B)$ of order $v$ such that every $T\in\binom{X}{t}$ is
contained in exactly $\lambda$ blocks of $\B$. To avoid triviality, we impose the following
restrictions on a $t$-$(v,k,\lambda)$ design $(X,\B)$:
\begin{enumerate}[(i)]
\item $t\geq 2$,
\item $t<k$,
\item $\B\not=\text{\O}$, and $\B\not=\binom{X}{k}$.
\end{enumerate}

For two designs, $(X,\mathcal{A})$ and $(Y,\mathcal{B})$,
  an \emph{isomorphism} of $(X,\mathcal{A})$ onto
  $(Y,\mathcal{B})$ is a bijection $\sigma:X\rightarrow Y$
  such that $\sigma(\mathcal{A})=\mathcal{B}$.
  An \emph{automorphism} of a design is an isomorphism of the design
  onto itself. The set of all automorphisms of a design $\mathcal{D}$
  forms a group under functional composition. This group is
  called the \emph{automorphism group} of $\mathcal{D}$ and is
  denoted by $\text{Aut}(\mathcal{D})$. A subgroup
  $H\leq \text{Aut}(\mathcal{D})$ is a \emph{group of automorphisms} 
  of $\mathcal{D}$.

Let $V$ be a set of cardinality $n$ and consider the induced
  action of the symmetric group $\mathcal{S}_n=\text{Sym}(V)$ on the set
  $X=\binom{V}{2}$. This defines an embedding of $\mathcal{S}_n$
  into $\mathcal{S}_{\binom{n}{2}}=\text{Sym}(X)$ with image group
  $\mathcal{S}_n^{[2]}$. By canonical extension, $\mathcal{S}_n^{[2]}$ 
  also acts on $\binom{X}{k}$.
  A $t$-$(v,k,\lambda)$ design $(X,\B)$ is \emph{graphical} if it has
  a group of automorphisms that is permutation isomorphic
  to $\mathcal{S}_n^{[2]}$ with $v=\binom{n}{2}$. In particular,
  $\mathcal{B}$ is then a union of orbits of $\mathcal{S}_n^{[2]}$
  on $\binom{X}{k}$.

  The term ``graphical design'' is motivated by the following alternative 
  perspective. Considering the complete graph $K_n$ with vertex set $V$,
  we may view $X$ as the edge set of $K_n$, in which case the orbits of 
  $\mathcal{S}_n^{[2]}$ on $\binom{X}{k}$ are in a one-to-one 
  correspondence with the isomorphism classes of spanning 
  $k$-edge subgraphs of $K_n$. Thus, we may view the block set 
  $\mathcal{B}$ of a graphical design as a set of spanning $k$-edge 
  subgraphs of $K_n$, closed under isomorphism of graphs, such that 
  every $t$-edge subgraph of $K_n$ is a subgraph of $\lambda$ graphs in 
  $\mathcal{B}$. Although the definition of a graphical design 
  does not explicitly assume this graphical structure, a required 
  group of automorphisms induces the structure (in a canonical
  manner for $n\neq 4$) because one of the orbits of $\mathcal{S}_n^{[2]}$ 
  corresponds to the line graph of $K_n$, from which one can recover
  the sets of edges having a vertex in common when $n\neq 4$.
  Two graphical designs 
  $(X,\mathcal{A})$ and $(Y,\mathcal{B})$, with individualized 
  required groups of automorphisms, $H$ and $K$, respectively, are 
  \emph{equivalent} if there exists an isomorphism 
  $\sigma$ of $(X,\mathcal{A})$ onto $(Y,\mathcal{B})$ 
  such that $\sigma H\sigma^{-1}=K$.

The first example of a graphical design has been attributed to R. M. Wilson
by Kramer and Mesner \cite{KramerMesner:1976}: \\

\noindent {\bf Example 1.1.}
{\em
A graphical \/ $3$-$(10,4,1)$ design is obtained by taking as blocks all
spanning $4$-edge subgraphs of $K_5$ isomorphic to one of
the following graphs: 
\begin{eqnarray*}
\begin{array}{ccccc}
\includegraphics[width=0.6in]{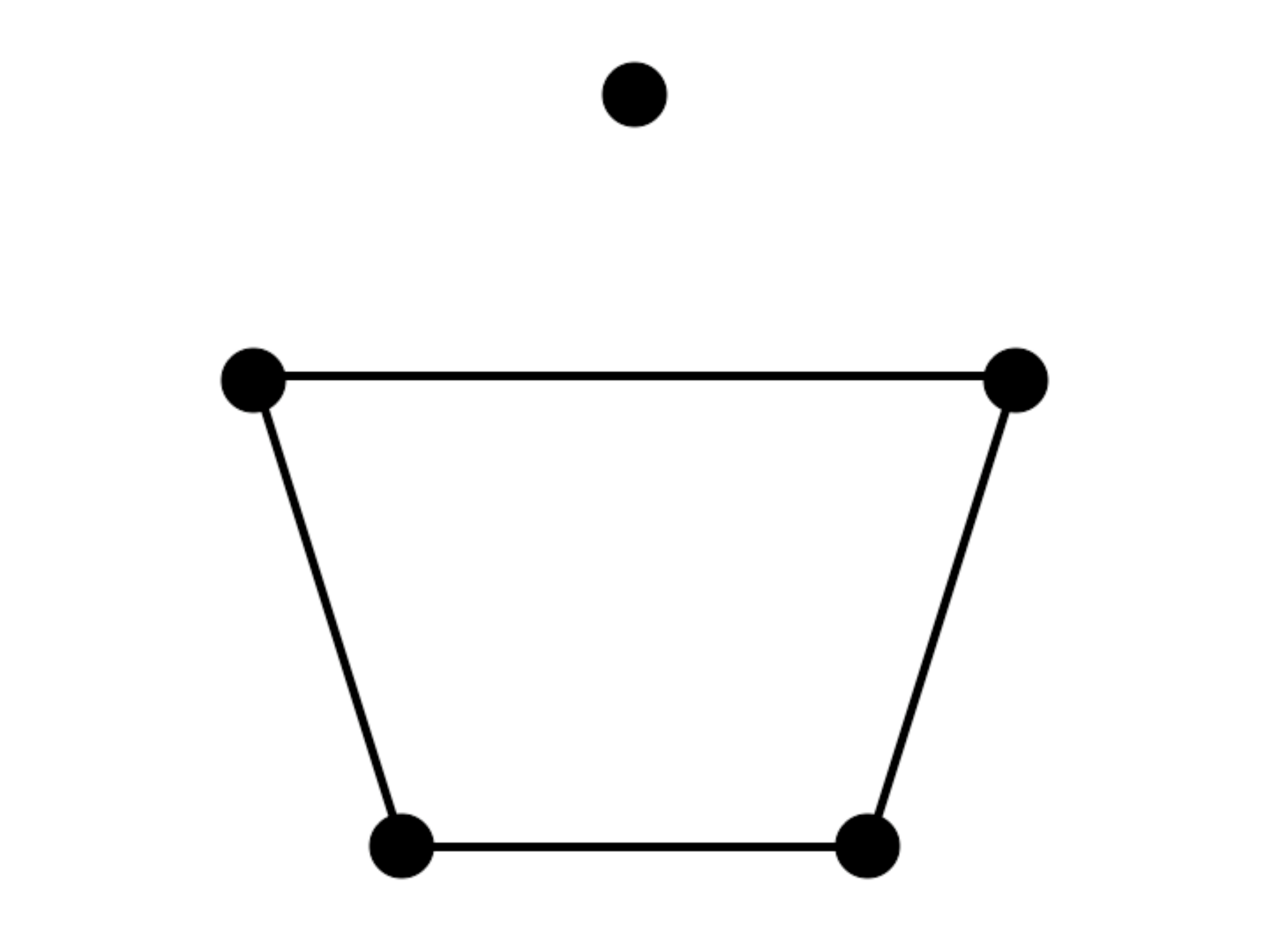} & &
\includegraphics[width=0.6in]{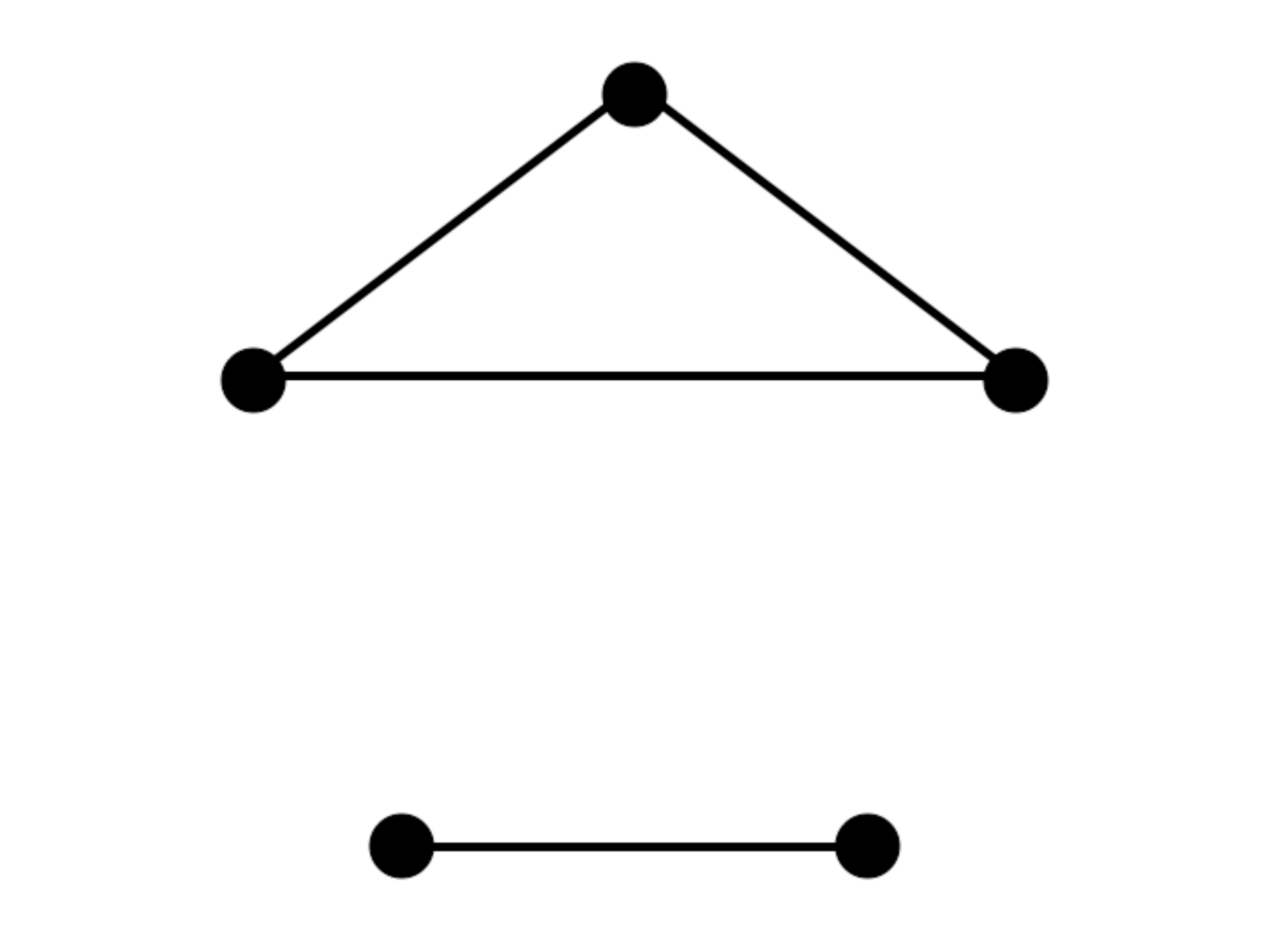} & &
\includegraphics[width=0.6in]{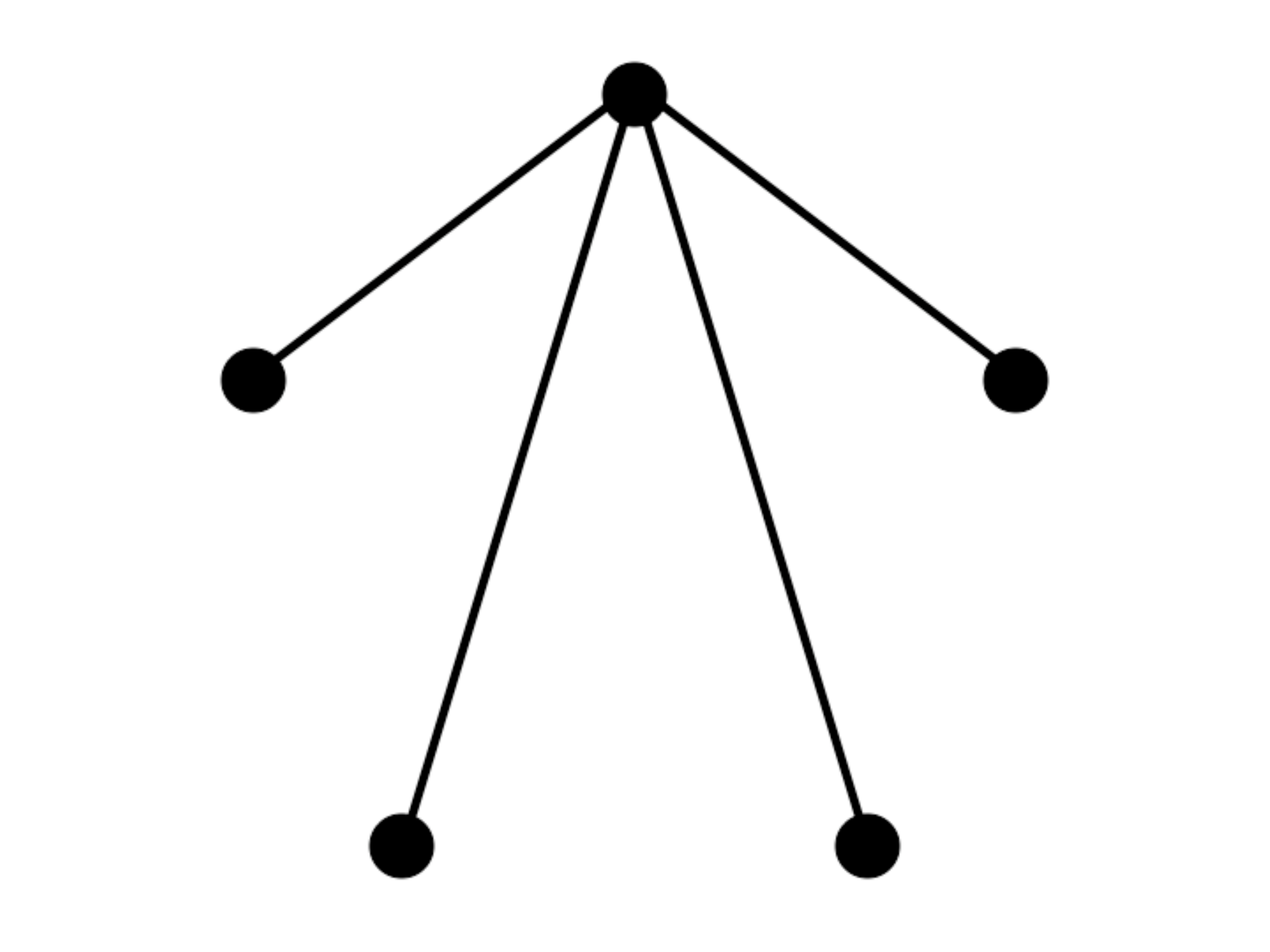}
\end{array}
\end{eqnarray*}
} \\

\noindent However, Betten {\em et al.}~\cite{Bettenetal:1999} have reported that already in 1970,
M. H. Klin has described graphical designs when he determined the overgroups
of $\SSS_n^{[2]}$. But Klin's result was unpublished, except for a short note that appeared in
a less well known journal \cite{Klin:1970}. Further examples of graphical designs were given
by Driessen \cite{Driessen:1978}. The first systematic approach to determining the existence
of graphical designs was undertaken by Chouinard {\em et al.} \cite{Chouinardetal:1983}, who
determined all graphical $t$-$(v,k,\lambda)$ designs with $\lambda=1$ and $\lambda=2$.
These results led Chouinard \cite{Chouinard:Private1989}
to make the following conjecture, which remains open. \\

\noindent {\bf Conjecture 1.2~(Chouinard).}
{\em
For any fixed $\lambda$, there exist only finitely many graphical $t$-$(v,k,\lambda)$ designs.}
\\

\noindent Partial progress on this conjecture has been obtained by Chouinard \cite{Chouinard:1996}.

Computers were brought to bear in the early nineties, which resulted in further progress in the
construction of graphical $t$-$(v,k,\lambda)$ designs. Kreher {\em et al.} \cite{Kreheretal:1990}
used the LLL algorithm to construct many examples of graphical $t$-$(v,k,\lambda)$ designs.
Chee \cite{Chee:1991,Chee:1992}
used symbolic computational methods to find all graphical 2-$(v,3,\lambda)$, 2-$(v,4,\lambda)$,
3-$(v,4,\lambda)$, and 4-$(v,5,\lambda)$ designs. Further sporadic examples were also obtained
by Kramer \cite{Kramer:1990} and Chee \cite{Chee:1993}. In the late nineties, more graphical
$t$-$(v,k,\lambda)$ designs were discovered by Betten {\em et al.} \cite{Bettenetal:1999} using
an improved implementation of the LLL algorithm. This is the state-of-the-art.
Despite that more than twenty years have passed since the introduction of graphical designs, only
a small finite number of them are known.
Let $\Lambda(t,k,v)$ denote the set of $\lambda\leq \frac{1}{2}\binom{v-t}{k-t}$
for which a graphical $t$-$(v,k,\lambda)$
design exists, and let $\Psi(t,k)=\{(v,\lambda):\lambda\in\Lambda(t,k,v)\}$.
The reason for restricting $\lambda\leq\frac{1}{2}\binom{v-t}{k-t}$ is to avoid duplication by
complementation, since if $(X,\B)$ is a (graphical) $t$-$(v,k,\lambda)$ design,
then its {\em complement},
$(X,\binom{X}{k}\setminus \B)$, is a (graphical) $t$-$(v,k,\binom{v-t}{k-t}-\lambda)$ design.
The parameters of all graphical designs known are given in Appendix A, where
Table \ref{complete} presents those sets $\Psi(t,k)$ which we have complete knowledge of,
and Table \ref{partial} lists known elements of some $\Psi(t,k)$ which we have yet to completely
determine. The authority for these tables are
\cite{Bettenetal:1999,Chee:1991,Chee:1992,Chee:1993,Kramer:1990,Kreheretal:1990}
(cf. \cite{CheeKreher:2006}).
In total, there are only
574 parameter sets for which we know there exist graphical designs.
Indeed, results in the literature are either on construction of sporadic examples, on nonexistence, or
on the finiteness of the number of graphical designs with certain parameters.

The purpose of this paper is to improve this state of knowledge by
determining completely the sets $\Psi(2,5)$ and $\Psi(3,5)$. With this result, the sets
$\Psi(t,k)$ are now completely known for $2\leq t<k\leq 5$. As a by-product, we give
more than 8000 new parameter sets for which there exists a graphical design,
substantially improving on the number of graphical
designs known thus far. Our results also correct some minor errors in \cite{Bettenetal:1999}.

\section{Kramer--Mesner Matrices and Outline of Approach}

Suppose we wish to construct a $t$-$(v,k,\lambda)$ design 
  $(X,\mathcal{B})$ with a group of automorphisms $\Gamma$.
  Then $\mathcal{B}$ is a union of orbits of $\Gamma$ on $\binom{X}{k}$.
  Let 
  $\mathcal{O}_1^{(t)},\mathcal{O}_2^{(t)},\ldots,\mathcal{O}_{N(t)}^{(t)}$
  and
  $\mathcal{O}_1^{(k)},\mathcal{O}_2^{(k)},\ldots,\mathcal{O}_{N(k)}^{(k)}$
  be the orbits of $\Gamma$ on $\binom{X}{t}$ and on $\binom{X}{k}$, 
  respectively. Define an $N(t)\times N(k)$ integer matrix 
  $W_{t,k}(X|\Gamma)$ by the rule that the $(i,j)$-entry is
  $|\{K\in\mathcal{O}_j^{(k)}:K\supseteq T\}|$, where $T\in\OO_i^{(t)}$
  can be chosen arbitrarily.
Such $W_{t,k}(X|\Gamma)$ matrices are called {\em Kramer--Mesner matrices}, after Kramer
and Mesner \cite{KramerMesner:1976} who observed the following.

\begin{theorem}[(Kramer and Mesner)]
\label{KM}
There exists a $t$-$(v,k,\lambda)$ design with a group of automorphisms $\Gamma$ if and only if
there exists a $\{0,1\}$-vector $\vu$ such that
\begin{eqnarray}
\label{KMeq}
W_{t,k}(X|\Gamma) \vu & = & \lambda (1,\ldots,1)^{\sf T}.
\end{eqnarray}
\end{theorem}

Based on Theorem \ref{KM},
our approach to determining $\Psi(2,5)$ and $\Psi(3,5)$ is to find all solutions to the
equation $W_{t,k}(X|\SSS_n^{[2]}) \vu = \lambda (1,\ldots,1)^{\sf T}$ for
$(t,k)=(2,5)$ and $(t,k)=(3,5)$.
More precisely, we perform the following steps:
\begin{enumerate}[(i)]
\item determine a bound $n_0$ so that no graphical $t$-$(v,k,\lambda)$ design
exists for $n\geq n_0$; and
\item enumerate all graphical $t$-$(\binom{n}{2},k,\lambda)$ designs for $n<n_0$ by
determining all solutions to $W_{t,k}(X|\SSS_n^{[2]}) \vu = \lambda (1,\ldots,1)^{\sf T}$.
\end{enumerate}
The first step is accomplished via a combinatorial analysis and the
second step is accomplished via computation.
It is not hard to see that
distinct $\{0,1\}$-vectors $\vu$
satisfying (\ref{KMeq}) give inequivalent graphical designs.

Betten \emph{et al.}~\cite{Bettenetal:1999} have computed the matrices
  $W_{2,5}(X|\mathcal{S}_n^{[2]})$ and
  $W_{3,5}(X|\mathcal{S}_n^{[2]})$.
  These take the forms given in
  Figs.~\ref{W25} and \ref{W35}, where $n^{\underline{k}}$ denotes 
  the \emph{falling factorial} $n(n-1)\cdots(n-k+1)$.
  Observe that the matrices have constant row sum
  $\binom{\binom{n}{2}-t}{k-t}$.
  A list of orbit representatives indexing the rows and columns of
  $W_{2,5}(X|\mathcal{S}_n^{[2]})$ and
  $W_{3,5}(X|\mathcal{S}_n^{[2]})$ 
  is given in Appendix B.

{\footnotesize
\begin{figure}[tb]
\renewcommand\arraystretch{1.5}
\begin{eqnarray*}
\begin{array}{cc}
\text{index}
\end{array} \\
\left[
\begin{array}{cc}
4(n-3) & 4 \\
\frac{5}{2}(n-3)^{\underline{3}} & 10(n-4)^{\underline{2}} \\
6(n-3)^{\underline{2}} & 16(n-4) \\
7(n-3)^{\underline{2}} & 12(n-4) \\
4(n-3)^{\underline{2}} & 4(n-4) \\
\frac{2}{3}(n-3)^{\underline{3}} & 4(n-4)^{\underline{2}} \\
\frac{1}{8}(n-3)^{\underline{4}} & \frac{7}{6}(n-4)^{\underline{3}} \\
\frac{2}{3}(n-3)^{\underline{4}} & 4(n-4)^{\underline{3}} \\
5(n-3)^{\underline{3}} & 20(n-4)^{\underline{2}} \\
\frac{3}{2}(n-3)^{\underline{3}} & 4(n-4)^{\underline{2}} \\
\frac{3}{2}(n-3)^{\underline{4}} & 14(n-4)^{\underline{3}} \\
4(n-3)^{\underline{3}} & 24(n-4)^{\underline{2}} \\
\frac{1}{4}(n-3)^{\underline{5}} & 4(n-4)^{\underline{4}} \\
6(n-3)^{\underline{2}} & 16(n-4) \\
(n-3)^{\underline{2}} & 4(n-4) \\
\frac{3}{2}(n-3)^{\underline{4}} & 14(n-4)^{\underline{3}} \\
5(n-3)^{\underline{3}} & 20(n-4)^{\underline{2}} \\
2(n-3)^{\underline{4}} & 12(n-4)^{\underline{3}} \\
\frac{7}{3}(n-3)^{\underline{3}} & 4(n-4)^{\underline{2}} \\
\frac{1}{48}(n-3)^{\underline{6}} & \frac{3}{4}(n-4)^{\underline{5}} \\
\frac{1}{4}(n-3)^{\underline{5}} & 4(n-4)^{\underline{4}} \\
\frac{1}{8}(n-3)^{\underline{5}} & \frac{7}{6}(n-4)^{\underline{4}} \\
\frac{1}{2}(n-3)^{\underline{3}} & 3(n-4)^{\underline{2}} \\
\frac{1}{4}(n-3)^{\underline{4}} & \frac{2}{3}(n-4)^{\underline{3}} \\
\frac{1}{6}(n-3)^{\underline{3}} & 0 \\
0 & \frac{1}{48}(n-4)^{\underline{6}}
\end{array}
\right] 
\begin{array}{c}
1 \\
2 \\
3 \\
4 \\
5 \\
6 \\
7 \\
8 \\
9 \\
10 \\
11 \\
12 \\
13 \\
14 \\
15 \\
16 \\
17 \\
18 \\
19 \\
20 \\
21 \\
22 \\
23 \\
24 \\
25 \\
26
\end{array}
\end{eqnarray*}
\caption{Transpose of the Kramer--Mesner matrix $W_{2,5}(X|\SSS_n^{[2]})$}
\label{W25}
\end{figure}
}

{\footnotesize
\begin{figure}[tb]
\renewcommand\arraystretch{1.5}
\begin{eqnarray*}
\begin{array}{c}
\text{index}
\end{array} \\
\left[
\begin{array}{ccccc}
3(n-3) & 0 & 3 & 3 & 0 \\
\frac{3}{2}(n-3)^{\underline{3}} & 5(n-5) & (n-4)^{\underline{2}} & \frac{3}{2}(n-4)^{\underline{2}} & 12 \\
3(n-3)^\2 & 8 & 4(n-4) & 3(n-4) & 0 \\
3(n-3)^\2 & 4 & 5(n-4) & 6(n-4) & 0 \\
\frac{3}{2}(n-3)^\2 & 1 & 2(n-4) & 6(n-4) & 0 \\
\frac{1}{2}(n-3)^\3 & 3(n-5) & 0 & 0 & 0 \\
\frac{1}{8}(n-3)^\4 & \frac{1}{2}(n-5)^\2 & 0 & 0 & 3(n-6) \\
0 & 3(n-5)^\2 & 0 & \frac{1}{2}(n-4)^\3 & 0 \\
0 & 12(n-5) & 3(n-4)^\2 & 3(n-4)^\2 & 0 \\
0 & 2(n-5) & (n-4)^\2 & \frac{3}{2}(n-4)^\2 & 0 \\
0 & 7(n-5)^\2 & \frac{1}{2}(n-4)^\3 & 0 & 24(n-6) \\
0 & 12(n-5) & 3(n-4)^\2 & 0 & 24 \\
0 & \frac{3}{2}(n-5)^\3 & 0 & 0 & 12(n-6)^\2 \\
0 & 6 & 6(n-4) & 3(n-4) & 0 \\
0 & 2 & n-4 & 0 & 0 \\
0 & 5(n-5)^\2 & (n-4)^\3 & 0 & 36(n-6) \\
0 & 8(n-5) & 4(n-4)^\2 & 3(n-4)^\2 & 24 \\
0 & 5(n-5)^\2 & (n-4)^\3 & \frac{3}{2}(n-4)^\3 & 24(n-6) \\
0 & 2(n-5) & (n-4)^\2 & 4(n-4)^\2 & 0 \\
0 & \frac{1}{8}(n-5)^\4 & 0 & 0 & \frac{7}{2}(n-6)^\3 \\
0 & (n-5)^\3 & \frac{1}{8}(n-4)^\4 & 0 & 15(n-6)^\2 \\
0 & \frac{1}{2}(n-5)^\3 & 0 & \frac{1}{8}(n-4)^\4 & 3(n-6)^\2 \\
0 & n-5 & \frac{1}{2}(n-4)^\2 & 0 & 6 \\
0 & \frac{1}{2}(n-5)^\2 & 0 & \frac{1}{2}(n-4)^\3 & 0 \\
0 & 0 & 0 & \frac{1}{2}(n-4)^\2 & 0 \\
0 & 0 & 0 & 0 & \frac{1}{8}(n-6)^\4 
\end{array}
\right]
\begin{array}{c}
1 \\
2 \\
3 \\
4 \\
5 \\
6 \\
7 \\
8 \\
9 \\
10 \\
11 \\
12 \\
13 \\
14 \\
15 \\
16 \\
17 \\
18 \\
19 \\
20 \\
21 \\
22 \\
23 \\
24 \\
25 \\
26
\end{array}
\end{eqnarray*}
\caption{Transpose of the Kramer--Mesner matrix $W_{3,5}(X|\SSS_n^{[2]})$}
\label{W35}
\end{figure}
}

\section{Upper Bounds for Existence}

Our subsequent proofs of the nonexistence of graphical
  designs for $n$ large enough in the cases $(t,k)=(2,5)$ and 
  $(t,k)=(3,5)$ are quantitative versions of the proof of a finiteness
  theorem of Betten \emph{et al.} \cite{Bettenetal:1999}.

  The orbit of a graph $G$ under the action of $\mathcal{S}_n^{[2]}$
  is denoted by $\text{Orb}(G)$.

\subsection{\boldmath Upper Bound for Existence of
Graphical $2$-$(v,5,\lambda)$ Designs}

We prove in this section that no graphical 2-$(\binom{n}{2},5,\lambda)$ design exists
if $n\geq 538$.

Let $(X,\B)$ be a graphical 2-$(\binom{n}{2},5,\lambda)$ design, where $n\geq 538$.
We may assume
without loss of generality that $\B\supseteq Orb(G_{26}^{(5)})$, since
otherwise we can consider the complement of the design. Let $\mu_i$ denote the sum
of all entries of degree $i$ (as a polynomial in $n$) in row two of $W_{2,5}(X|\SSS_n^{[2]})$.
Then we have
$\mu_6  =  \frac{1}{48}(n-4)^\6$,
$\mu_5 =  \frac{3}{4}(n-4)^\5$,
$\mu_4  =  \frac{55}{6}(n-4)^\4$,
$\mu_3  =  \frac{275}{6}(n-4)^\3$, and
$\mu_2  =  89(n-4)^\2$.
Define the integers $\lambda_6=\mu_6$ and $\lambda_i = \lambda_{i+1}+\mu_i$ for $i=2,3,4,5$.
By considering the number of blocks
in $Orb(G_{26}^{(5)})$ containing $G_2^{(2)}$, 
we see that
\begin{eqnarray}
\label{lambda21}
\lambda & \geq & \lambda_6.
\end{eqnarray}

\begin{lemma}
$\B\supseteq Orb(G_{20}^{(5)})$.
\end{lemma}

\begin{proof}
Suppose that $\B\not\supseteq Orb(G_{20}^{(5)})$. Then by considering the number of blocks
in $\B$ containing $G_1^{(2)}$, we have
\begin{eqnarray*}
\lambda & \leq &  \binom{\binom{n}{2}-2}{3}-\frac{1}{48}(n-3)^\6.
\end{eqnarray*}
The above inequality, together with inequality (\ref{lambda21}), implies
\begin{eqnarray*}
\lambda_6 & \leq & \binom{\binom{n}{2}-2}{3}-\frac{1}{48}(n-3)^\6,
\end{eqnarray*}
giving
\begin{eqnarray*}
n^6-69n^5+1085n^4-8435n^3+36642n^2-84664n+80832 & \leq & 0,
\end{eqnarray*}
which is impossible for $n\geq 51$.
\end{proof} \\

So $\B\supseteq\bigcup_{i\in\{20,26\}} Orb(G_i^{(5)})$ and by considering the number
of blocks in $\B$ containing $G_2^{(2)}$, we now have
\begin{eqnarray}
\label{lambda22}
\lambda & \geq & \lambda_5.
\end{eqnarray}

\begin{lemma}
$\B\supseteq \bigcup_{i\in\{13,21,22\}} Orb(G_i^{(5)})$.
\end{lemma}

\begin{proof}
Suppose that $\B$ contains at most two of the orbits $Orb(G_i^{(5)})$, $i\in\{13$, $21$, $22\}$.
Then by considering the number of blocks in $\B$ containing $G_1^{(2)}$, we have
\begin{eqnarray*}
\lambda & \leq & \binom{\binom{n}{2}-2}{3} - \frac{1}{8}(n-3)^\5.
\end{eqnarray*}
The above inequality, together with inequality (\ref{lambda22}), implies
\begin{eqnarray*}
\lambda_5 & \leq & \binom{\binom{n}{2}-2}{3} - \frac{1}{8}(n-3)^\5,
\end{eqnarray*}
giving
\begin{eqnarray*}
3n^5 - 295n^4 + 4475n^3 - 28541n^2 + 85198n - 98184 & \leq & 0,
\end{eqnarray*}
which is impossible for $n\geq 82$.
\end{proof} \\

So $\B\supseteq\bigcup_{i\in\{13,20,21,22,26\}} Orb(G_i^{(5)})$ and by considering the number
of blocks in $\B$ containing $G_2^{(2)}$, we now have
\begin{eqnarray}
\label{lambda23}
\lambda & \geq & \lambda_4.
\end{eqnarray}

\begin{lemma}
$\B\supseteq \bigcup_{i\in\{7,8,11,16,18,24\}} Orb(G_i^{(5)})$.
\end{lemma}

\begin{proof}
Suppose that $\B$ contains at most five of the orbits $Orb(G_i^{(5)})$,
$i\in\{7$, $8$, $11$, $16$, $18$, $24\}$.
Then by considering the number of blocks in $\B$ containing $G_1^{(2)}$, we have
\begin{eqnarray*}
\lambda & \leq & \binom{\binom{n}{2}-2}{3} - \frac{1}{8}(n-3)^\4.
\end{eqnarray*}
The above inequality, together with inequality (\ref{lambda23}), implies
\begin{eqnarray*}
\lambda_4 & \leq & \binom{\binom{n}{2}-2}{3} - \frac{1}{8}(n-3)^\4,
\end{eqnarray*}
giving
\begin{eqnarray*}
3n^4 - 1154n^3 + 14721n^2 - 64450n + 95256 & \leq & 0,
\end{eqnarray*}
which is impossible for $n\geq 372$.
\end{proof} \\

So $\B\supseteq\bigcup_{i\in\{7,8,11,13,16,18,20,21,22,24,26\}} Orb(G_i^{(5)})$ and by considering the number
of blocks in $\B$ containing $G_2^{(2)}$, we now have
\begin{eqnarray}
\label{lambda24}
\lambda & \geq & \lambda_3.
\end{eqnarray}

\begin{lemma}
$\B\supseteq \bigcup_{i\in\{2,6,9,10,12,17,19,23,25\}} Orb(G_i^{(5)})$.
\end{lemma}

\begin{proof}
Suppose that $\B$ contains at most eight of the orbits $Orb(G_i^{(5)})$, $i\in\{2$, 6, 9,
10, 12, 17, 19, 23, $25\}$.
Then by considering the number of blocks in $\B$ containing $G_1^{(2)}$, we have
\begin{eqnarray*}
\lambda & \leq & \binom{\binom{n}{2}-2}{3} - \frac{1}{6}(n-3)^\3.
\end{eqnarray*}
The above inequality, together with inequality (\ref{lambda24}), implies 
\begin{eqnarray*}
\lambda_3 & \leq & \binom{\binom{n}{2}-2}{3} - \frac{1}{6}(n-3)^\3,
\end{eqnarray*}
giving
\begin{eqnarray*}
n^3 - 546n^2 + 4541n - 9516 & \leq & 0,
\end{eqnarray*}
which is impossible for $n\geq 538$.
\end{proof} \\

So $\B\supseteq\bigcup_{i\in\{2,6,7,8,9,10,11,12,13,16,17,18,19,20,21,22,23,24,25,26\}} Orb(G_i^{(5)})$ and by considering the number
of blocks in $\B$ containing $G_2^{(2)}$, we now have
\begin{eqnarray}
\label{lambda25}
\lambda & \geq & \lambda_2.
\end{eqnarray}

\begin{lemma}
$\B\supseteq \bigcup_{i\in\{3,4,5,14,15\}} Orb(G_i^{(5)})$.
\end{lemma}

\begin{proof}
Suppose that $\B$ contains at most four of the orbits $Orb(G_i^{(5)})$, $i\in\{3$, 4, 5, 14, $15\}$.
Then by considering the number of blocks in $\B$ containing $G_1^{(2)}$, we have
\begin{eqnarray*}
\lambda & \leq & \binom{\binom{n}{2}-2}{3} - (n-3)^\2.
\end{eqnarray*}
The above inequality, together with inequality (\ref{lambda25}), implies
\begin{eqnarray*}
\lambda_2 & \leq & \binom{\binom{n}{2}-2}{3} - (n-3)^\2,
\end{eqnarray*}
giving
\begin{eqnarray*}
n^2-59n+216 & \leq & 0,
\end{eqnarray*}
which is impossible for $n\geq 56$.
\end{proof} \\

So $\B\supseteq\bigcup_{i\in\{2,3,4,5,6,7,8,9,10,11,12,13,14,15,16,17,18,19,20,21,22,23,24,25,26\}} Orb(G_i^{(5)})=\binom{X}{5}\setminus Orb(G_1^{(5)})$. If $\B\not\supseteq Orb(G_1^{(5)})$,
then $(X,\B)$ cannot be a $2$-$(\binom{n}{2},5,\lambda)$
design unless $4(n-3)=4$, which is impossible for $n\geq 5$. So
$\B\supseteq Orb(G_1^{(5)})$ and hence $\B=\binom{X}{k}$, which is excluded from the
definition of a design to avoid triviality.

We summarize the above results as:

\begin{theorem}
No graphical $2$-$(\binom{n}{2},5,\lambda)$ design exists if
$n\geq 538$.
\end{theorem}

\subsection{\boldmath Upper Bound for Existence of
Graphical $3$-$(v,5,\lambda)$ Designs}

We prove in this section that no graphical 3-$(\binom{n}{2},5,\lambda)$ design exists
if $n\geq 34$.

Let $(X,\B)$ be a graphical 3-$(\binom{n}{2},5,\lambda)$ design, where $n\geq 34$.
We may assume without loss of generality that $\B\supseteq Orb(G_7^{(5)})$, since
otherwise we can consider the complement of the design. By considering the number of blocks
in $Orb(G_7^{(5)})$ containing $G_1^{(3)}$, we see that
\begin{eqnarray}
\label{L1}
\lambda & \geq & \frac{1}{8}(n-3)^\4.
\end{eqnarray}

\begin{lemma}
$\B\supseteq\bigcup_{i\in\{20,21,22,26\}} Orb(G_i^{(5)})$.
\end{lemma}

\begin{proof}
Suppose that $\B\not\supseteq Orb(G_{20}^{(5)})$, Then by considering the number of blocks in
$\B$ containing $G_2^{(3)}$, we have
\begin{eqnarray*}
\lambda & \leq & \binom{\binom{n}{2}-3}{2} - \frac{1}{8}(n-5)^\4.
\end{eqnarray*}
The above inequality, together with inequality (\ref{L1}), implies 
\begin{eqnarray*}
 \frac{1}{8}(n-3)^\4 & \leq &  \binom{\binom{n}{2}-3}{2} - \frac{1}{8}(n-5)^\4,
\end{eqnarray*}
giving
\begin{eqnarray*}
(n-4)(n^3-38n^2+231n-498) & \leq & 0,
\end{eqnarray*}
which is impossible for $n\geq 32$.

To show that $\B\supseteq Orb(G_i^{(5)})$ for $i\in\{21,22,26\}$, mimic the proof above.
\end{proof} \\

It follows that $\B\supseteq\bigcup_{i\in\{7,20,21,22,26\}}Orb(G_i^{(5)})$.
Let $\A=\binom{X}{5}\setminus\B$ and
consider the 3-$(\binom{n}{2},5,\lambda')$ design $(X,\A)$. By
considering the number
of blocks in $\A$ containing $G_5^{(3)}$, we see that
\begin{eqnarray}
\label{L'}
\lambda' & \leq & 12(n-6)^\2+84(n-6)+66.
\end{eqnarray}

\begin{lemma}
$\A\not\supseteq Orb(G_i^{(5)})$ for $i\in\{2,6,11,13,16,18,24\}$.
\end{lemma}

\begin{proof}
Suppose that $\A\supseteq Orb(G_2^{(5)})$. Then by considering the number of blocks in
$Orb(G_2^{(5)})$ containing $G_1^{(3)}$, we have
\begin{eqnarray*}
\lambda' & \geq & \frac{3}{2}(n-3)^\3.
\end{eqnarray*}
The above inequality, together with inequality (\ref{L'}), implies 
\begin{eqnarray*}
 \frac{3}{2}(n-3)^\3 & \leq & 12(n-6)^\2+84(n-6)+66,
\end{eqnarray*}
giving
\begin{eqnarray*}
n^3 - 20n^2 + 95n - 120 & \leq & 0,
\end{eqnarray*}
which is impossible for $n\geq 14$.

To show that $A\not\supseteq Orb(G_i^{(5)})$ for $i\in\{6,11,13,16,18,24\}$, mimic
the proof above.
\end{proof} \\

It follows that $\B\supseteq\bigcup_{i\in\{2,6,7,11,13,16,18,20,21,22,24,26\}} Orb(G_i^{(5)})$.
By considering the number of blocks in $\A$ containing $G_5^{(3)}$, we now have
\begin{eqnarray}
\label{L''}
\lambda' & \leq & 54.
\end{eqnarray}

\begin{lemma}
$\A\not\supseteq Orb(G_i^{(5)})$ for $i\in\{1,3,4,5,8,9,10,12,14,15,17,19,23,25\}$.
\end{lemma}

\begin{proof}
Suppose that $\A\supseteq Orb(G_1^{(5)})$. Then by considering the number of blocks
in $Orb(G_1^{(5)})$ containing $G_1^{(3)}$, we have
\begin{eqnarray*}
\lambda' & \geq & 3(n-3).
\end{eqnarray*}
The above inequality, together with inequality (\ref{L''}), implies 
\begin{eqnarray*}
3(n-3) & \leq & 54,
\end{eqnarray*}
which is impossible for $n\geq 22$.

To show that $\A\not\supseteq Orb(G_i^{(5)})$ for $i\in\{3,4,5,8,9,10,12,14,15,17,19,23,25\}$,
mimic the proof above.
\end{proof} \\

We can now conclude that $\B\supseteq \binom{X}{5}$, which is excluded from the definition
of a design to avoid triviality. We summarize the above results as:

\begin{theorem}
No graphical $3$-$(\binom{n}{2},5,\lambda)$ design exists if
$n\geq 34$.
\end{theorem}

\section{Computation for Existence}

The symbolic computation approach of Chee \cite{Chee:1991} can, in theory, be
  used to find all graphical $t$-$(v,k,\lambda)$ designs for given
  $t$ and $k$, without the need to establish upper bounds for existence,
  such as in the previous section. However, in practice, the method
  becomes infeasible when $k$ becomes large. Already for $k=5$ we
  would have to solve up to 33 million systems of simultaneous 
  Diophantine equations of degree up to six. 
  Fortunately, using the upper bounds from the previous section,
  a straightforward exhaustive search suffices.
  In both of the cases $(t,k)=(2,5)$ and $(t,k)=(3,5)$,
  there are 26 possible orbits of 5-edge graphs, implying 
  that we can easily enumerate all the $2^{26}=67108864$ candidate 
  designs, represented as $\{0,1\}$-vectors $\mathsf{u}$, 
  and filter out those candidates that do not constitute a solution
  to the system 
  
  \begin{eqnarray*}
  W_{t,5}(X|\mathcal{S}_n^{[2]})\mathsf{u}=\lambda(1,\ldots,1)^\top,
  \quad
  \lambda\leq\frac{1}{2}\binom{\binom{n}{2}-t}{5-t}.
  \end{eqnarray*}
  
\noindent In particular, this system needs to be considered
  in the two cases $t=2$ and $t=3$ for all $n\leq 537$ and $n\leq 39$, 
  respectively. Both authors of this paper independently carried out 
  this computation with the following identical results.
    
 \subsection{\boldmath Existence of Graphical $2$-$(v,5,\lambda)$ Designs}
 
 Our computations show that there are no graphical 2-$(\binom{n}{2},5,\lambda)$
 designs for $40\leq n\leq 537$. For $n\leq 39$, the number of inequivalent
 graphical 2-$(\binom{n}{2},5,\lambda)$ designs is fairly large, and 
 for reasons of space, it is infeasible to give
 a complete listing within this paper. A complete catalogue of the designs 
can be found on the first author's website at
\begin{center}
 $\langle$ \url{http://www1.spms.ntu.edu.sg/~ymchee/graphical.htm} $\rangle$.
 \end{center}
 We record this result as:
 
 \begin{theorem}
 There are $8619$ elements in $\Psi(2,5)$ and
there exist $271360$ inequivalent
 graphical $2$-$(\binom{n}{2},5,\lambda)$ designs.
No graphical $2$-$(\binom{n}{2},5,\lambda)$ design exists
 if $n\geq 40$.
 \end{theorem}
 
\subsection{\boldmath Existence of Graphical $3$-$(v,5,\lambda)$ Designs}

Our computations show that there are no graphical 3-$(\binom{n}{2},5,\lambda)$
designs for $10\leq n\leq 33$. For $n\leq 9$, a complete listing of 
all inequivalent graphical 3-$(\binom{n}{2},5,\lambda)$
designs found is presented below.

\begin{longtable}{| c | c | c |}
\hline
\multicolumn{3}{| c |}{All elements of $\Psi(3,5)$ and inequivalent solutions} \\
\hline
$(v,\lambda)$ & $\{0,1\}$-vectors $\vu^{\sf T}$ giving & Number of  \\
& inequivalent solutions & inequivalent solutions \\
\hline
$(15,30)$ & 10010100110000001000001000 & 1 \\
\hline
$(21,3)$ & 00000010000000100000000010 & 1 \\
\hline
$(21,30)$ & 00001100001001000000001100 & 1 \\
\hline
$(21,33)$ & 00001110001001100000001110 & 1 \\
\hline
$(21,39)$ & 00010010100000010000000010 & 3 \\
                   & 01000011010000101000000000 & \\
                   & 01000011010001100100000000 & \\
\hline
$(21,48)$ & 10010000001100100010000010 & 2 \\
                   & 10100000000100100110000010 & \\
\hline
$(21,69)$ & 00011110101001010000001110 & 5 \\
                   & 00101110100000011000001110 & \\
                   & 00101110100001010100001110 & \\
                   & 00101111000101010010001110 & \\
                   & 01001111011001101000001100 & \\
\hline
$(21,75)$ & 01010011110000011000000000 & 2 \\
                   & 01010011110001010100000000 & \\
\hline
$(28,30)$ & 00000100000011000000001110 & 1 \\
\hline
$(28,150)$ & 00110101010100010110001000 & 4 \\
                     & 00110101011000011010001000 & \\
                     & 11001010100111100100110110 & \\
                     & 11001010101011101000110110 & \\
\hline
$(36,180)$ & 00101010011001000101001000 & 1 \\
\hline
$(36,198)$ & 11000000011001110010001110 & 1 \\
\hline
$(36,258)$ & 10101111000110111011000110 & 3 \\
                     & 10110100011100100110100010 & \\
                     & 10110100101010101010100010 & \\
\hline
\end{longtable}

\noindent We record this result as:

 \begin{theorem}
 There are $13$ elements in $\Psi(3,5)$ and
there exist $26$ inequivalent
 graphical $3$-$(\binom{n}{2},5,\lambda)$ designs.
No graphical $3$-$(\binom{n}{2},5,\lambda)$ design exists
 if $n\geq 10$.
 \end{theorem}
 
\section{Conclusion}

In this paper, we determined completely the sets $\Psi(2,5)$ and $\Psi(3,5)$, and found
more than 270000 inequivalent graphical designs, and more than 8000 new parameter sets
for which there exists a graphical design.

We remark that our computation revealed two minor errors in \cite{Bettenetal:1999};
in fact,
\begin{enumerate}[(i)]
\item there is only one graphical 2-$(21,5,\lambda)$ design for 
$\lambda=52$ and $\lambda=84$; and
\item there exist only two inequivalent (and hence at most two nonisomorphic)
graphical 3-$(21,5,75)$ designs.
\end{enumerate}

A natural question is whether the techniques in this paper could be developed further to
determine $\Psi(t,k)$ for higher $k$, in particular for $k=6$. The method for establishing
upper bounds for existence is certainly applicable, but the main hurdle is the
search for solutions to $W_{t,k}(X|\SSS_n^{[2]})\vu=\lambda(1,\ldots,1)^{\sf T}$.
There are 68 nonisomorphic $6$-edge graphs, so the na{\"\i}ve search space has size $2^{68}$.
More sophisticated search techniques must be employed in this case.

\appendix{A}
\section*{All Known Graphical \lowercase{$t$}-Designs}
\label{current}

{\normalsize
\begin{longtable}{| r | r | l | c |}
\caption[Complete knowledge of $\Psi(t,k)$]{Complete knowledge of $\Psi(t,k)$}
\label{complete} \\
\hline
$t$ & $k$ & All elements of $\Psi(t,k)$ & $|\Psi(t,k)|$ \\
\hline
2 & 3 &
$\begin{array}{ccccc}
(10,4) & (15,1) & (28,6) & (28,10) & (55,25)
\end{array}$ & 5 \\
\hline

2 & 4 &
$\begin{array}{ccccc}
(10,2) & (10,4) & (10,8) & (10,10) & (10,12) \\
(15,6) & (15,24) & (15,30 & (15,36) & (21,6) \\
(21,12) & (21,18) & (21,36) & (21,42) & (21,45) \\
(21,48) & (21,51) & (21,54) & (21,57) & (21,60) \\
(21,63) & (21,66) & (21,69) & (21,72) & (21,75) \\
(21,78) & (21,81) & (21,84) & (28,5) & (28,55) \\
(28,80) & (28,85) & (28,95) & (29,110) & (28,120) \\
(28,125) & (28,135) & (28,150) & (36,15) & (36,90) \\
(36,111) & (36,120) & (36,135) & (36,165) & (36,210) \\
(36,231) & (36,240) & (36,255) & (36,276) & (45,63) \\
(45,105) & (45,252) & (45,357) & (45,378) & (45,420) \\
(55,168) & (55,336) & (55,504) & (78,630) & (78,1080) \\
(78,1350) & (91,836) & (91,1430) & (91,1496) & (105,1320) \\
(105,1326) & (105,1650) & (105,1656) & (105,1782) & (105,1788) \\
(105,1980) & (105,1986) & (105,2112) & (105,2118) & (105,2442) \\
(105,2448) & (153,4935) & (153,5025) & (253,14535) &  
\end{array}$ & 79 \\
\hline

3 & 4 &
$\begin{array}{c}
(10,1) \\
\end{array}$ & 1 \\
\hline

4 & 5 &
$\begin{array}{c}
- \\
\end{array}$ & 0 \\
\hline

5 & 6 &
$\begin{array}{c}
- \\
\end{array}$ & 0 \\
\hline
\end{longtable}

\begin{longtable}{|r|r|l|c|}
\caption[Partial knowledge of $\Psi(t,k)$]{Partial knowledge of $\Psi(t,k)$}
\label{partial} \\
\hline
$t$ & $k$ & Known elements of $\Psi(t,k)$ & $|\Psi(t,k)|\geq $ \\
\hline
2 & 5 &
$\begin{array}{ccccc}
(10,16) & (10,20) & (21,7) & (21,12) & (21,19) \\
(21,22) & (21,34) & (21,35) & (21,47) & (21,50) \\
(21,52) & (21,55) & (21,57) & (21,60) & (21,62) \\
(21,64) & (21,67) & (21,69) & (21,70) & (21,72) \\
(21,77) & (21,79) & (21,82) & (21,84) & (21,89) \\
(21,94) & (21,95) & (21,100) & (21,120) & (28,60) \\
(28,100) & (28,140) & (28,160) & (28,200) & (28,240) \\
(28,260) & (28,300) & (28,340) & (28,360) & (36,60) \\
(36,80) & (36,140) & (36,164) & (36,180) & (36,224) \\
(36,240) & (36,244) & (36,480) & (36,720) & 
\end{array}$ & 98 \\
& & 
$(15,\lambda)$: $16\leq\lambda\leq 142$, $\lambda\equiv 0,2,4,$ or $6\pmod{10}$,
$\lambda\not=20,50$ & \\
\hline

2 & 6 &
$\begin{array}{ccccc}
(21,13) & (21,30) & (21,38) & (21,45) & (21,48) \\
(21,50) & (21,51) & (21,55) & (21,58) & (21,60) \\
(21,61) & (21,63) & (21,68) & (21,70) & (28,25) \\
(28,40) & (28,50) & (28,65) & (28,70) & (28,80) \\
(28,90) & (28,100) & (36,20) & (36,45) & (36,120) \\
(36,240) & (36,540) & (36,720) & (36,1080) & (36,2160) \\
(36,4320) & & & & 
\end{array}$ & 78 \\
& & 
$(15,\lambda)$: $10\leq\lambda\leq 355$, $\lambda\equiv 0$ or $10\pmod{15}$ & \\
\hline

2 & 7 &
$\begin{array}{ccccc}
(15,3) & (15,24) & (15,27) & (15,30) & (15,33) \\
(15,36) & (15,39) & (21,42) & (21,63) & (21,78) \\
(21,84) & (21,105) & (28,16) & (28,140) & (28,156) \\
(28,182) & (28,198) & (36,210) & (36,246) & (36,336) \\
(36,372) & (36,420) & (36,456) & (36,462) & (36,546) 
\end{array}$ & 224 \\
& & 
$(15,\lambda)$: $48\leq\lambda\leq 642$, $\lambda\equiv 0\pmod{3}$ & \\
\hline

2 & 8 &
$\begin{array}{ccccc}
(21,84) & (21,168) & (21,336) & (21,672) & (28,70) \\
(28,210) & & & & \\
\end{array}$ & 6 \\
\hline

2 & 9 &
$\begin{array}{ccccc}
(21,12) & (21,54) & (21,72) & (21,108) & (21,216) \\
(21,432) & (21,864) & (28,40) & (28,160) & (28,320) \\
(28,480) & (28,640) & (28,960) & (28,1920) & (28,3840)
\end{array}$ & 15 \\
\hline

3 & 5 &
$\begin{array}{ccccc}
(15,30) & (21,3) & (21,30) & (21,33) & (21,39) \\
(21,48) & (21,69) & (21,75) & (28,30) & (28,150) \\
(36,180) & (36,270) & & & 
\end{array}$ & 12 \\
\hline

3 & 6 &
$\begin{array}{ccccc}
(15,100) & (21,68) & (21,100) & (21,108) & (21,128) \\
(21,136) & (21,140) & (21,148) & (21,156) & (21,160) \\
(21,168) & (21,176) & (21,180) & (21,188) & (21,196) \\
(21,200) & (28,80) & (28,120) & (28,180) & (28,220) \\
(28,240) & (28,260) & & & 
\end{array}$ & 22 \\
\hline

3 & 7 &
$\begin{array}{ccccc}
(15,60) & (15,75) & (15,90) & (15,135) & (15,150) \\
(15,165) & (15,180) & (15,225) & (15,240) & (21,105) \\
(21,120) & (21,210) & (21,225) & (21,315) & (28,210) \\
(28,225) & (28,240) & (28,275) & & \\
\end{array}$ & 18 \\
\hline

3 & 8 &
$\begin{array}{ccccc}
(21,168) & (21,252) & (21,336) & (21,420) & (28,168) \\
(28,378) & (28,672) & & & 
\end{array}$ & 7 \\
\hline

3 & 9 &
$\begin{array}{ccccc}
(28,280) & & & & 
\end{array}$ & 1 \\
\hline

4 & 6 &
$\begin{array}{ccccc}
(28,132) & & & & \\
\end{array}$ & 1 \\
\hline

4 & 7 &
$\begin{array}{ccccc}
(15,60)
\end{array}$ & 1 \\
\hline

5 & 7 &
$\begin{array}{ccccc}
(28,93) & (36,165) & & &
\end{array}$ & 2 \\
\hline

5 & 8 &
$\begin{array}{ccccc}
(28,756) & (28,791) & (28,840) & (28, 875) & 
\end{array}$ & 4 \\
\hline

\end{longtable}
}

\appendix{B}
\section*{Orbit Representatives}
\label{orbitreps}

A list of orbit representatives for $t$-edge graphs, for $t=2$, $t=3$ and $t=5$, is given below.
Note that isolated vertices are not shown in our drawings. The orbit representative
indexing row $i$ of $W_{t,5}(X|\SSS_n^{[2]})$ is the graph $G_i^{(t)}$, $t\in\{2,3\}$,
and the orbit representative indexing column $j$ of $W_{t,5}(X|\SSS_n^{[2]})$
is the graph $G_j^{(5)}$.

\addtocounter{table}{2}
{\normalsize
\begin{longtable}{|c|c|}
\caption[Orbit representatives of 2-edge graphs]{Orbit representatives of 2-edge graphs} \\
\hline
$G_1^{(2)}$ & $G_2^{(2)}$ \\
\hline
\includegraphics[width=0.6in]{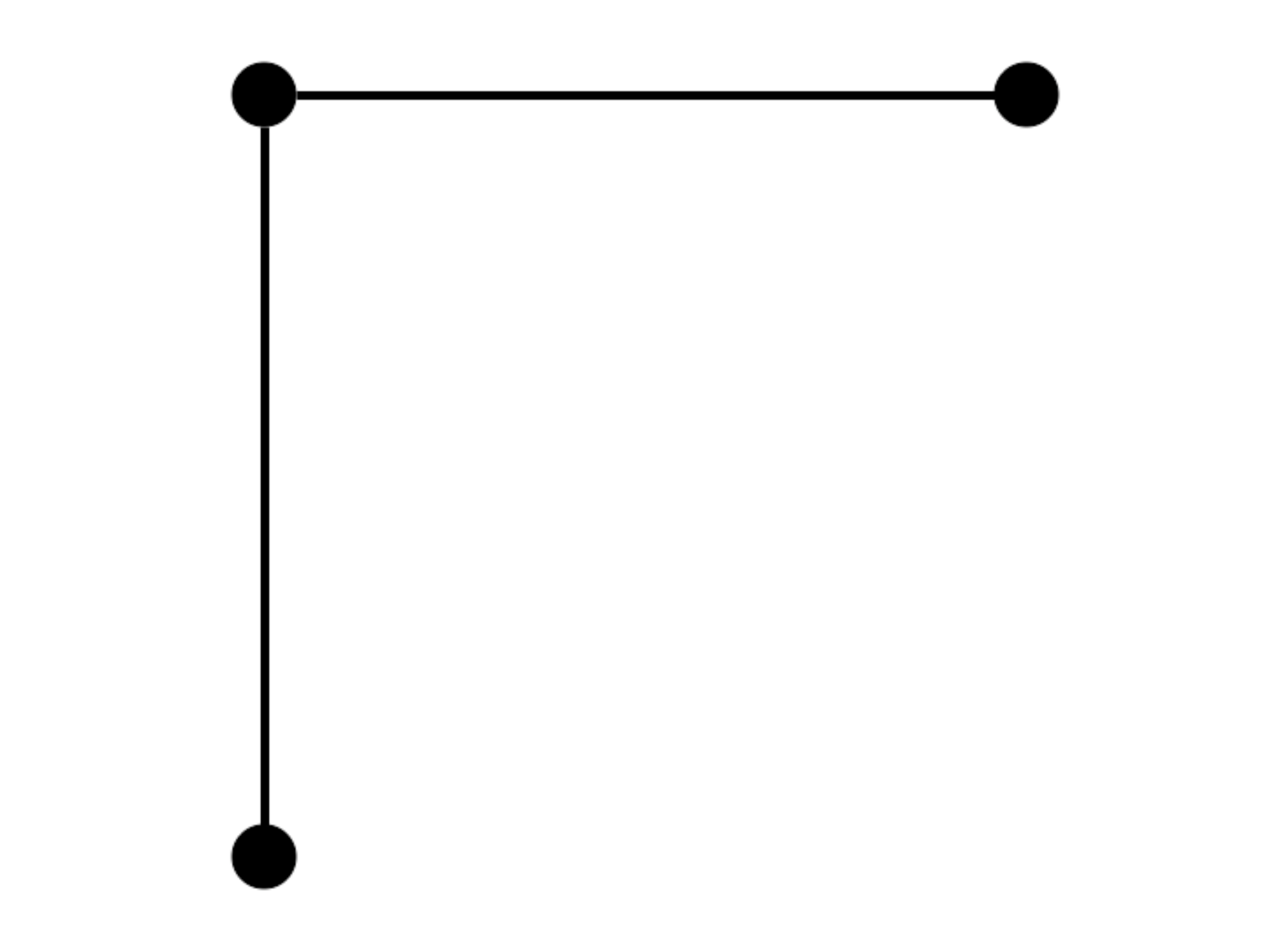} &
\includegraphics[width=0.6in]{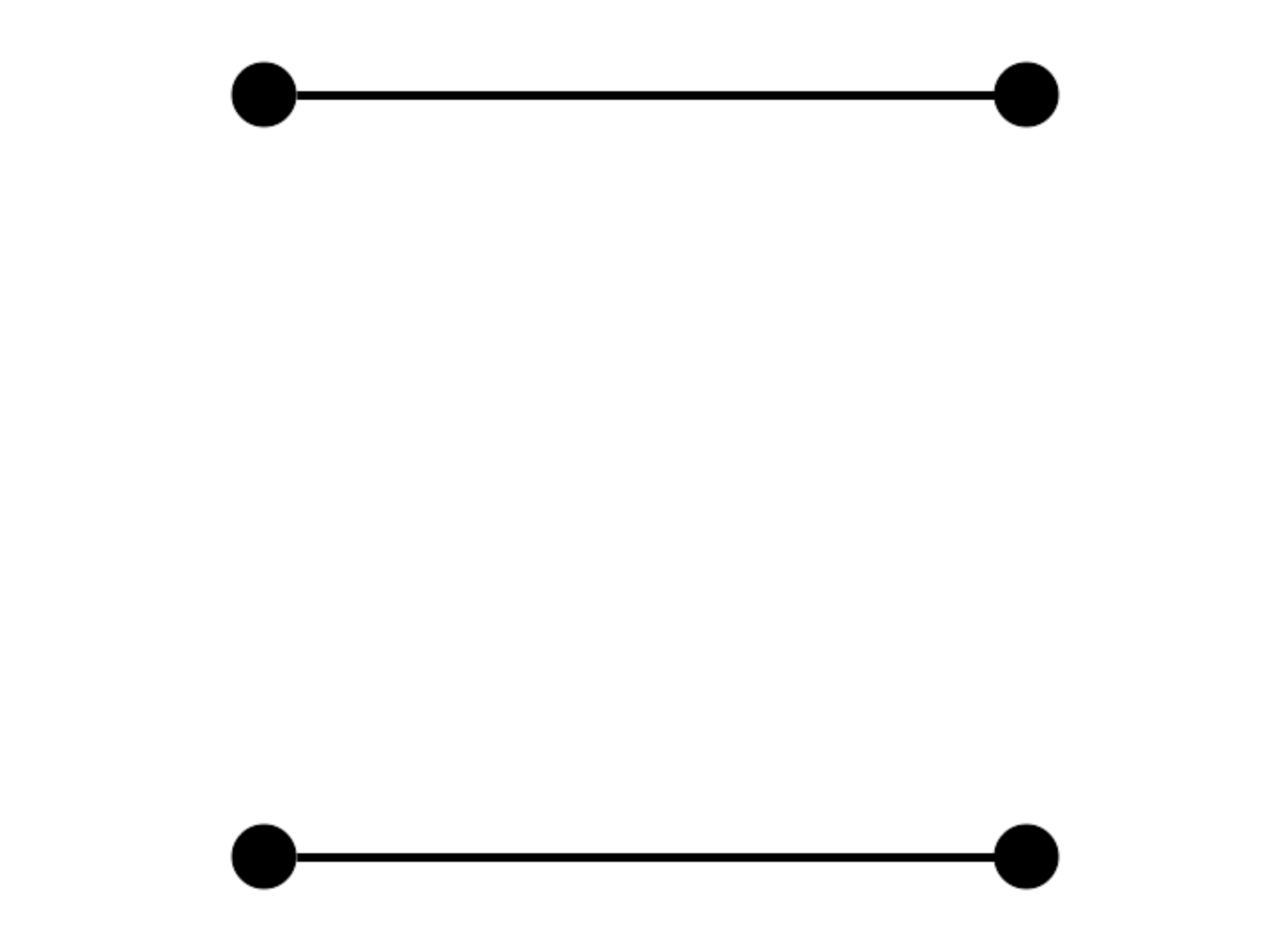} \\
\hline
\end{longtable}

\begin{longtable}{|c|c|c|c|c|}
\caption[Orbit representatives of 3-edge graphs]{Orbit representatives of 3-edge graphs} \\
\hline
$G_1^{(3)}$ & $G_2^{(3)}$ & $G_3^{(3)}$ & $G_4^{(3)}$ & $G_5^{(3)}$ \\
\hline
\includegraphics[width=0.6in]{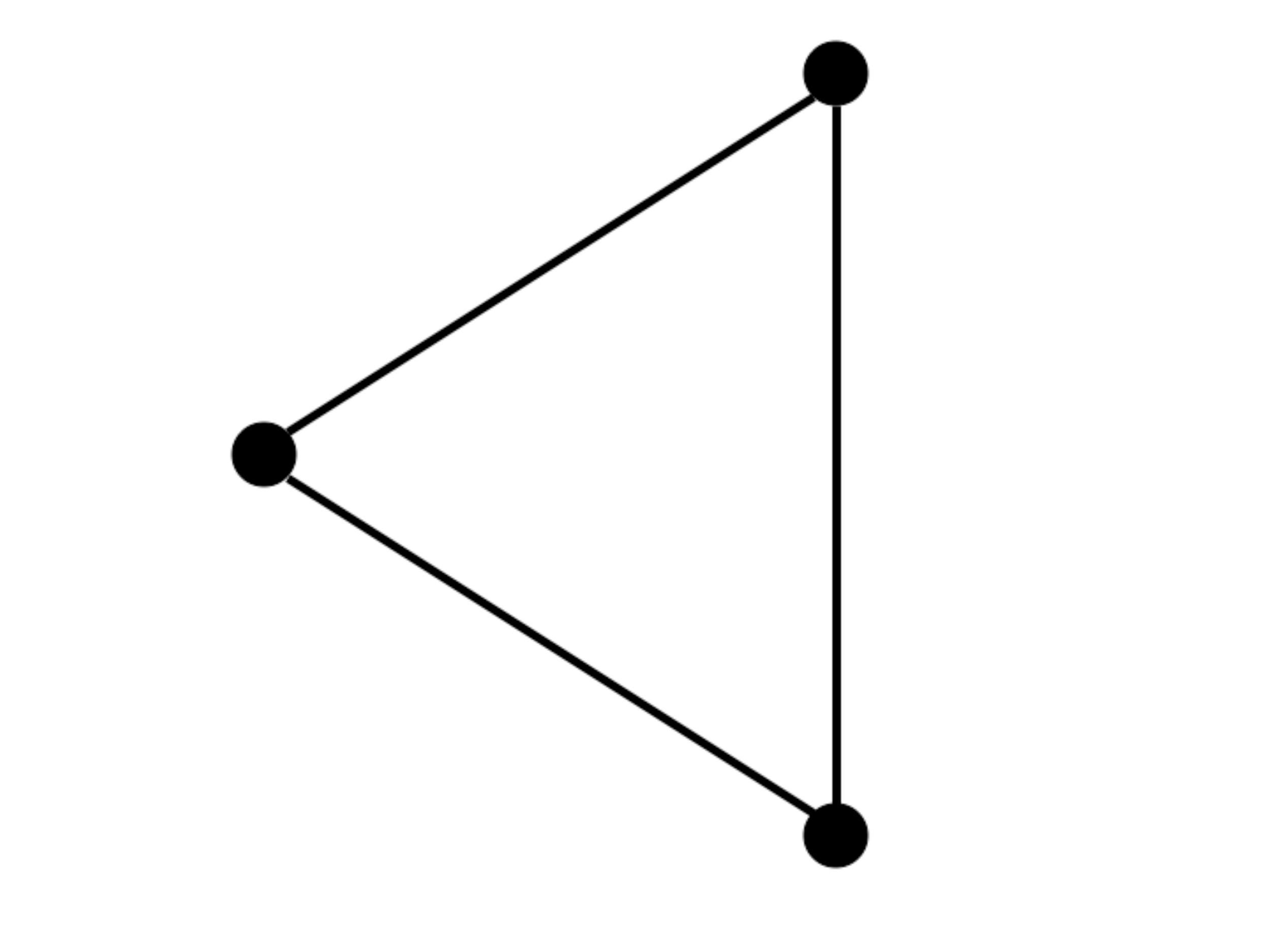} &
\includegraphics[width=0.6in]{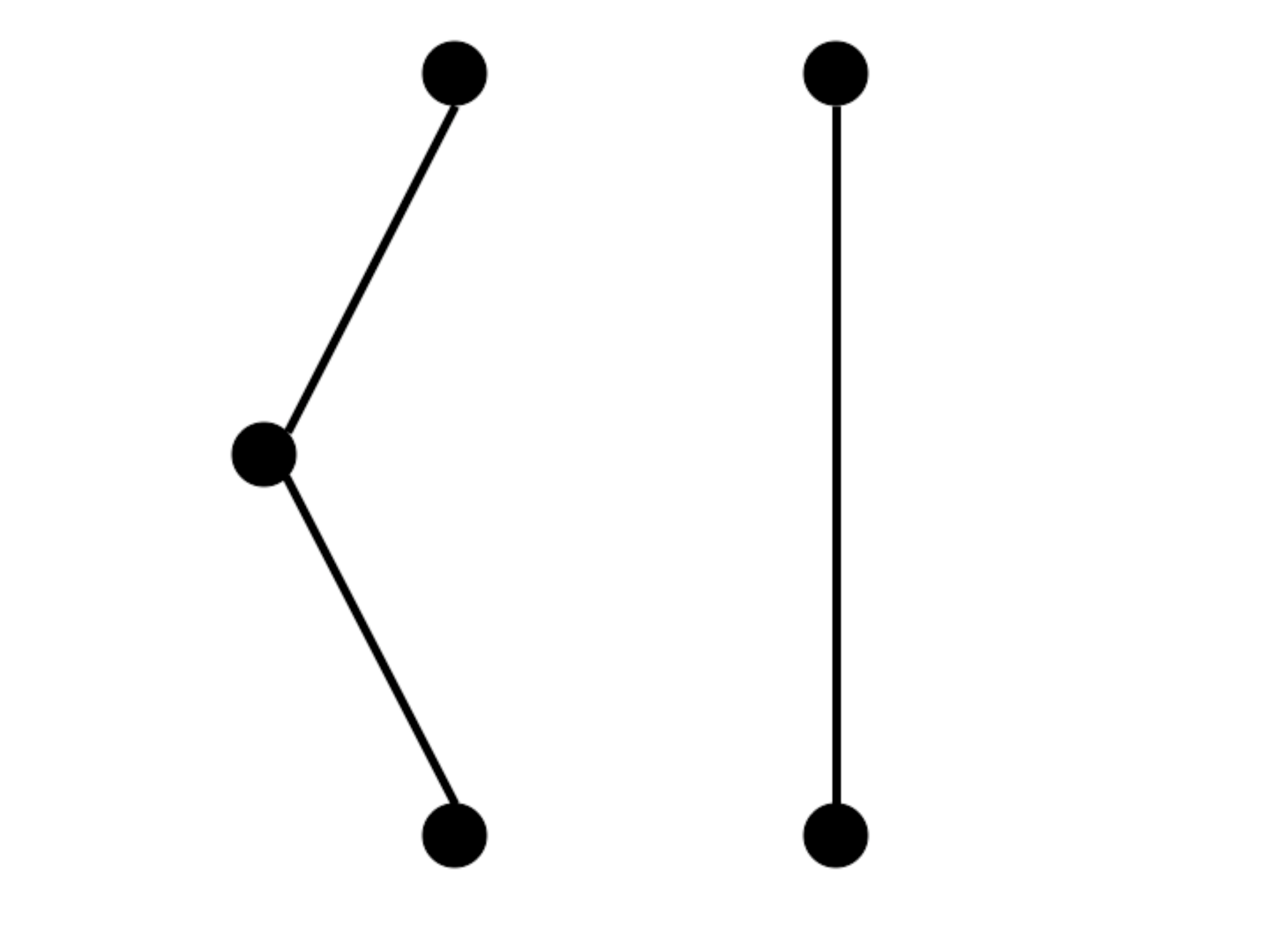} &
\includegraphics[width=0.6in]{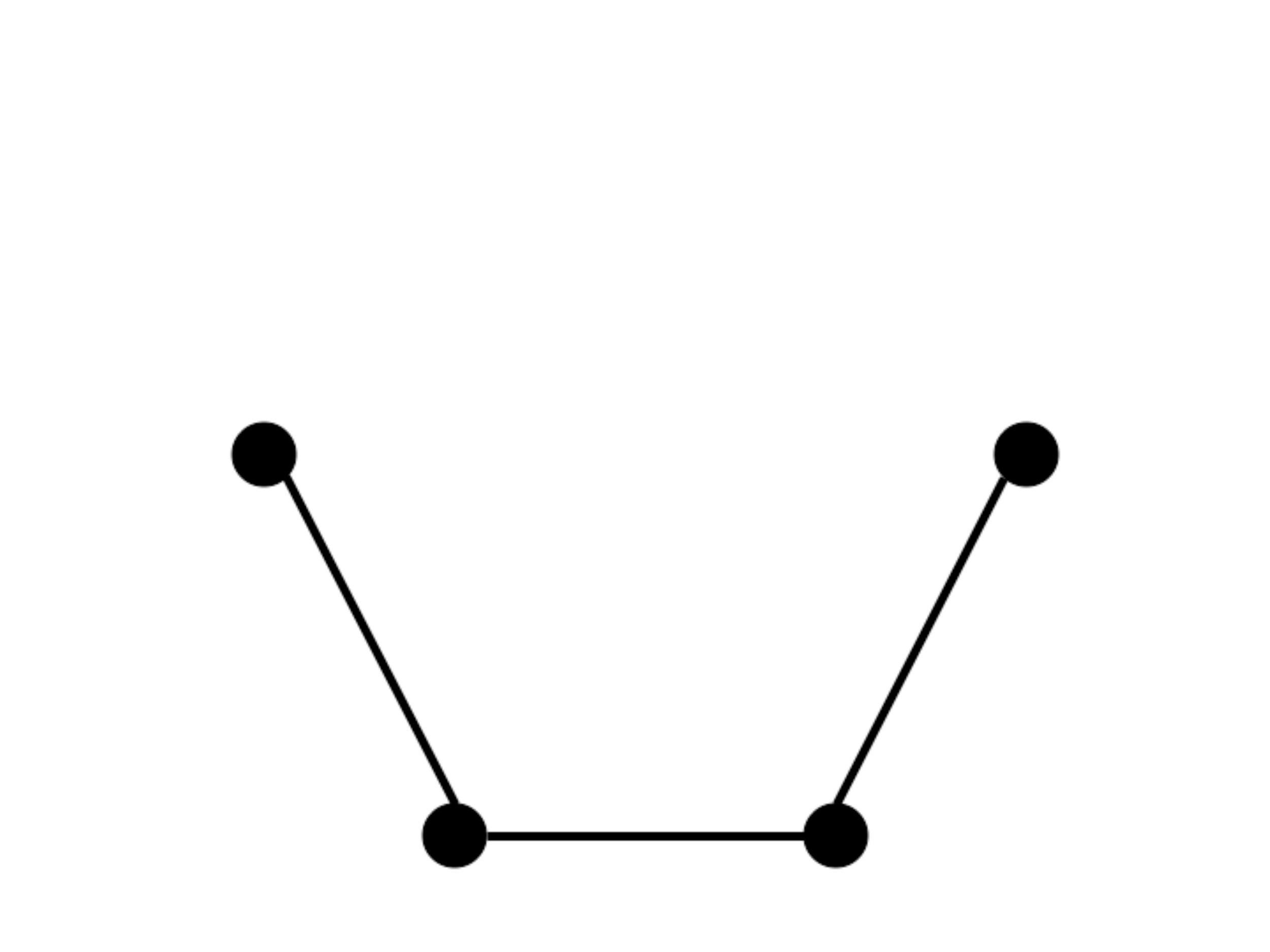} &
\includegraphics[width=0.6in]{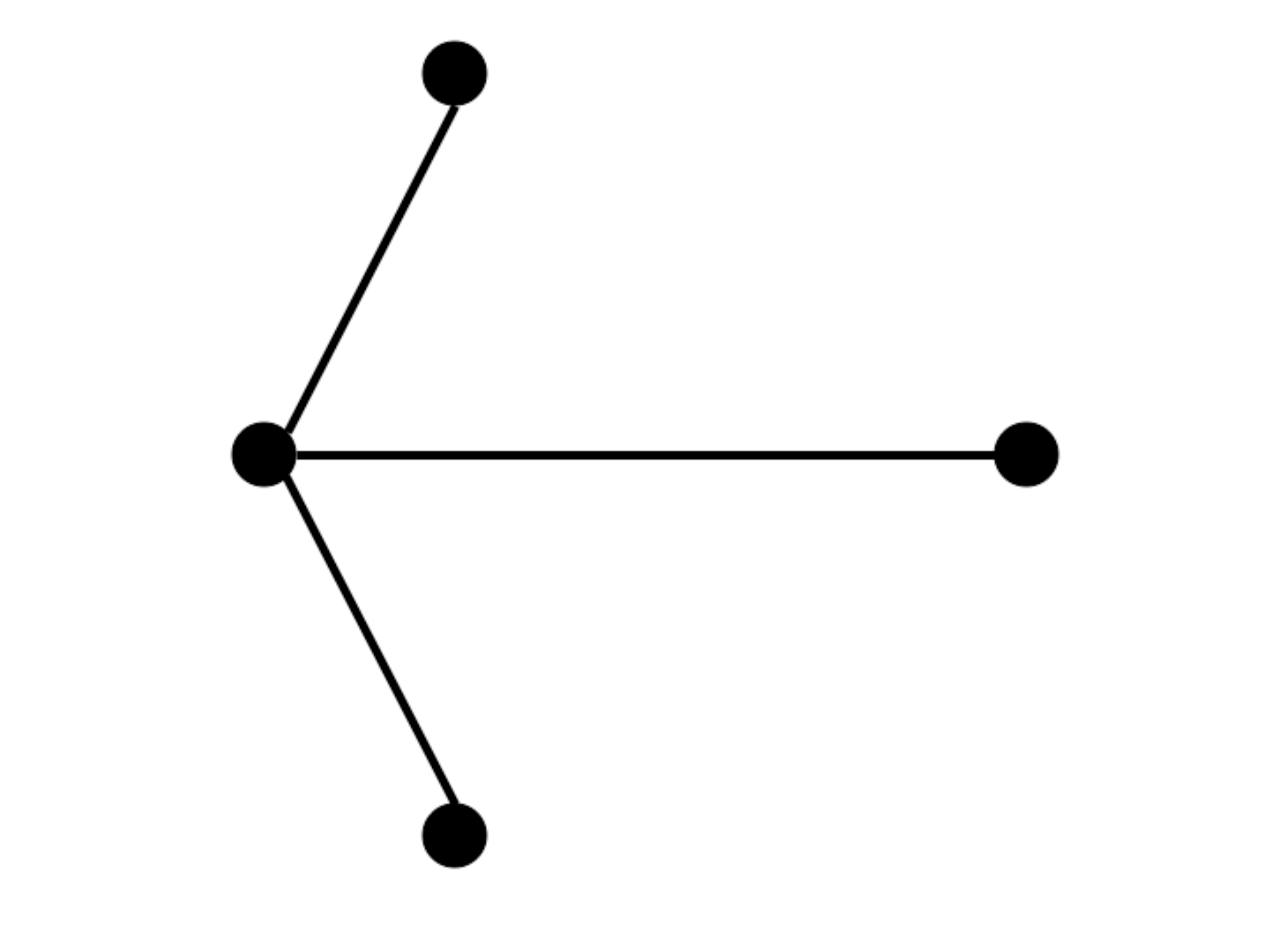} &
\includegraphics[width=0.6in]{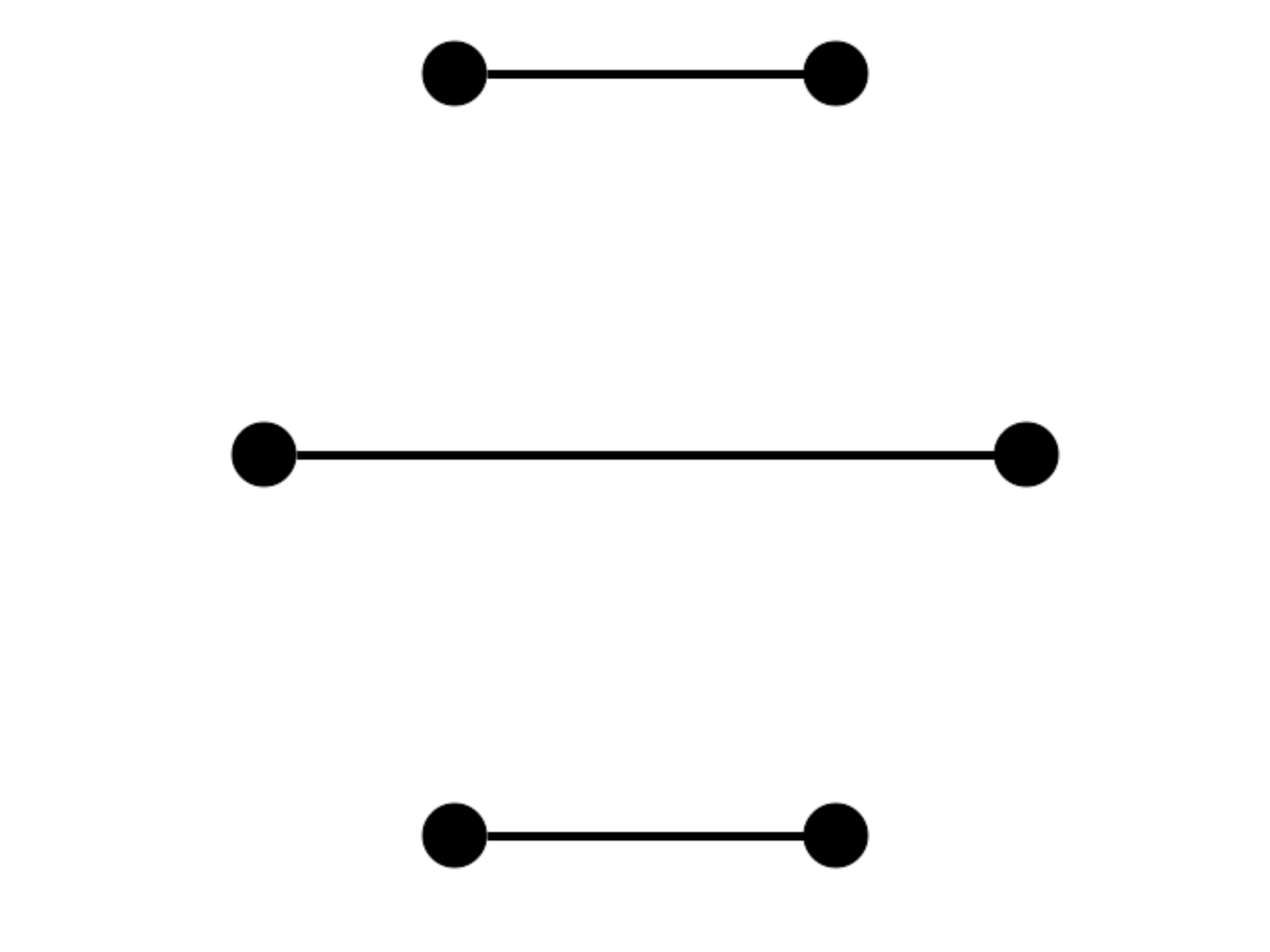} \\
\hline
\end{longtable}

\begin{longtable}{|c|c|c|c|c|c|}
\caption[Orbit representatives of 5-edge graphs]{Orbit representatives of 5-edge graphs} \\
\hline
$G_1^{(5)}$ & $G_2^{(5)}$ & $G_3^{(5)}$ & $G_4^{(5)}$ & $G_5^{(5)}$ & $G_6^{(5)}$ \\
\hline
\includegraphics[width=0.6in]{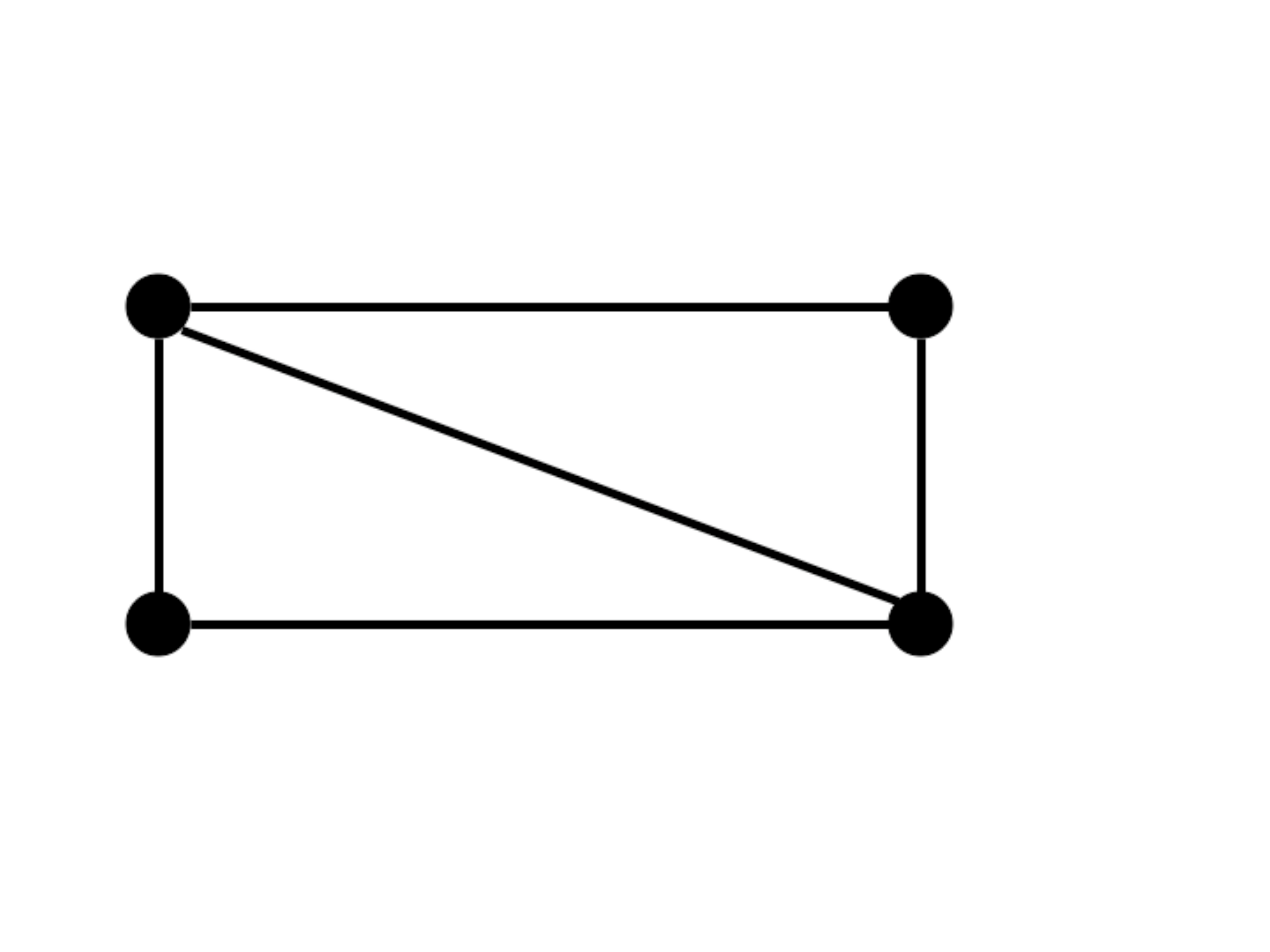} &
\includegraphics[width=0.6in]{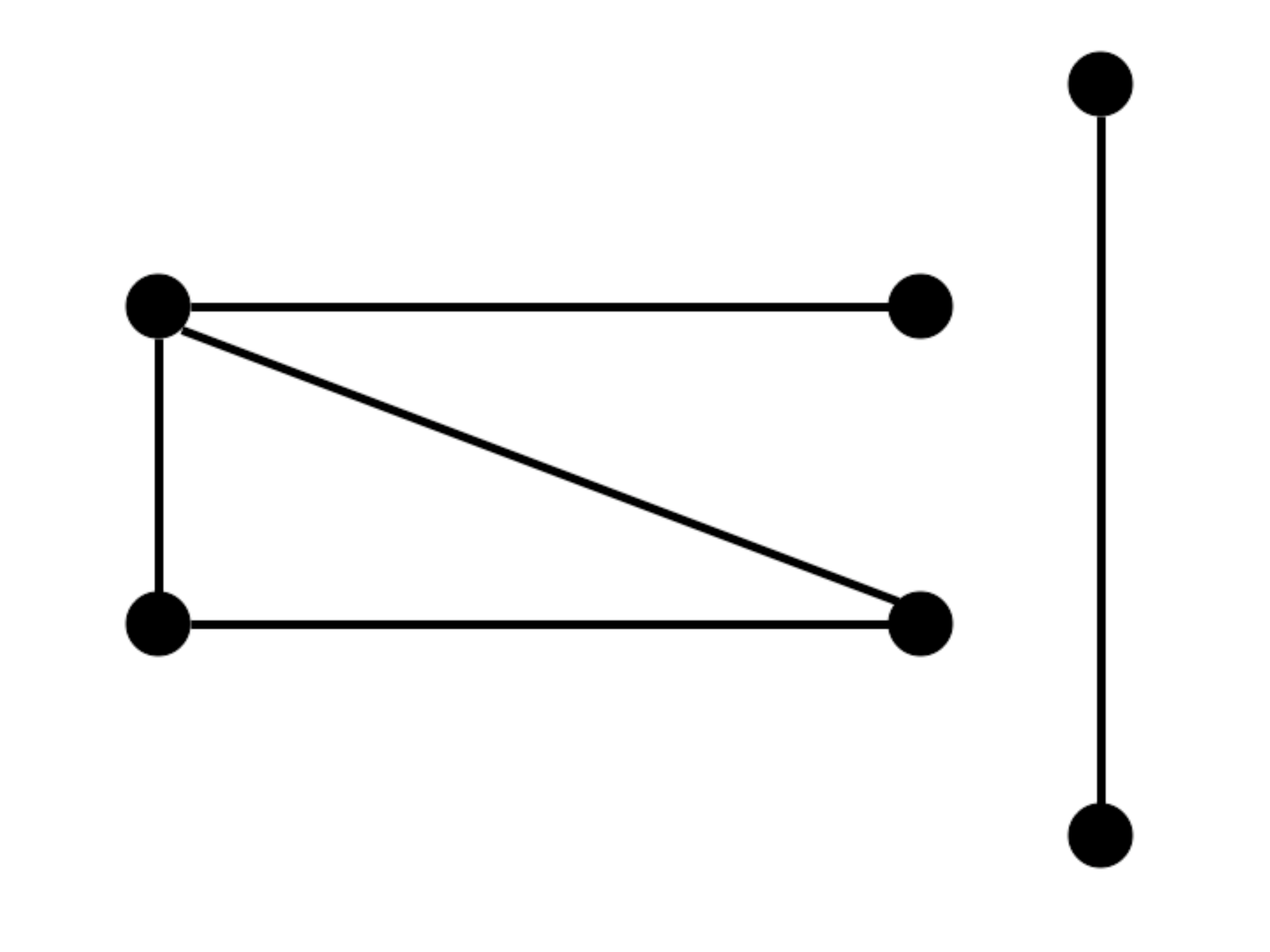} &
\includegraphics[width=0.6in]{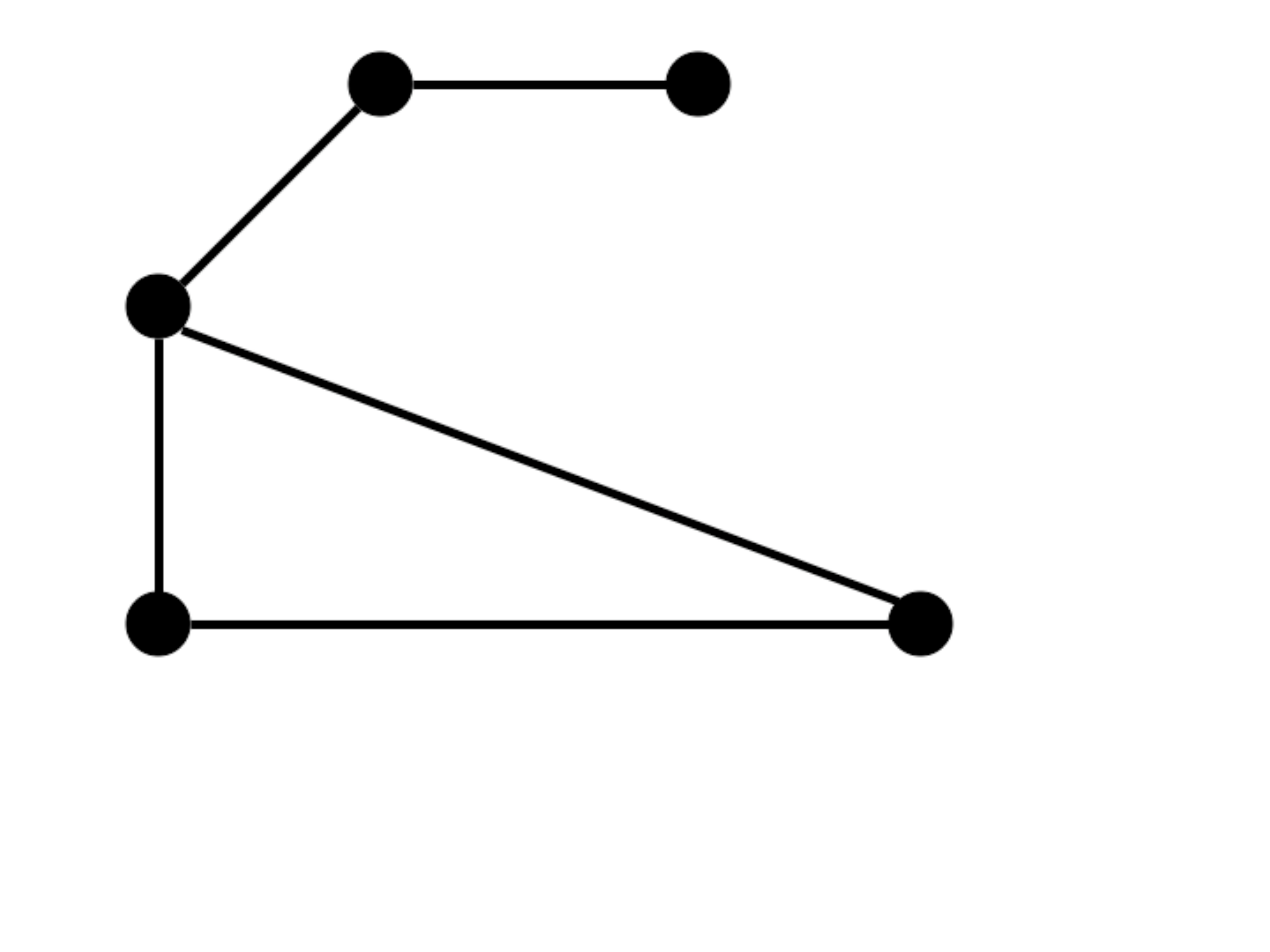} &
\includegraphics[width=0.6in]{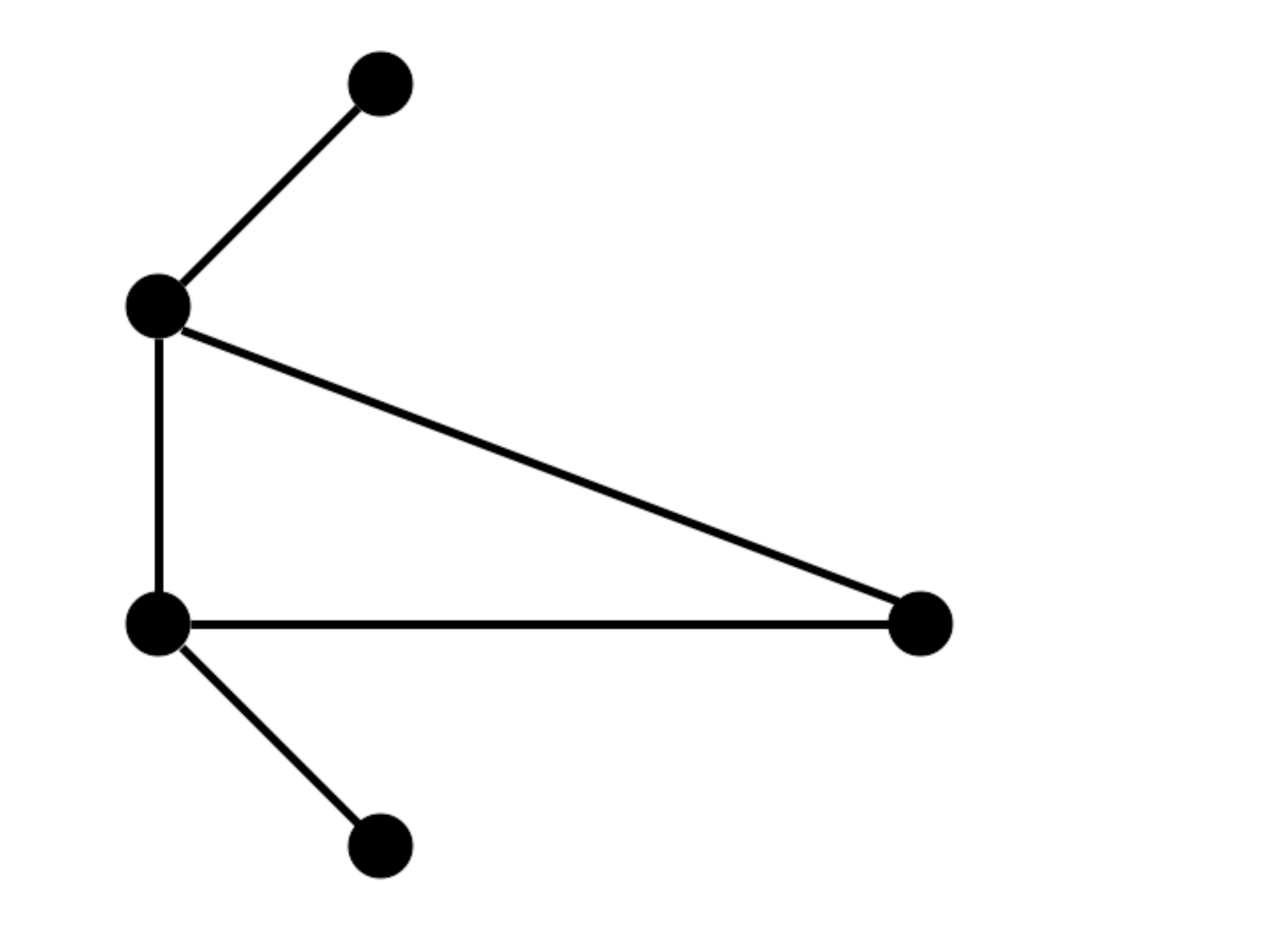} &
\includegraphics[width=0.6in]{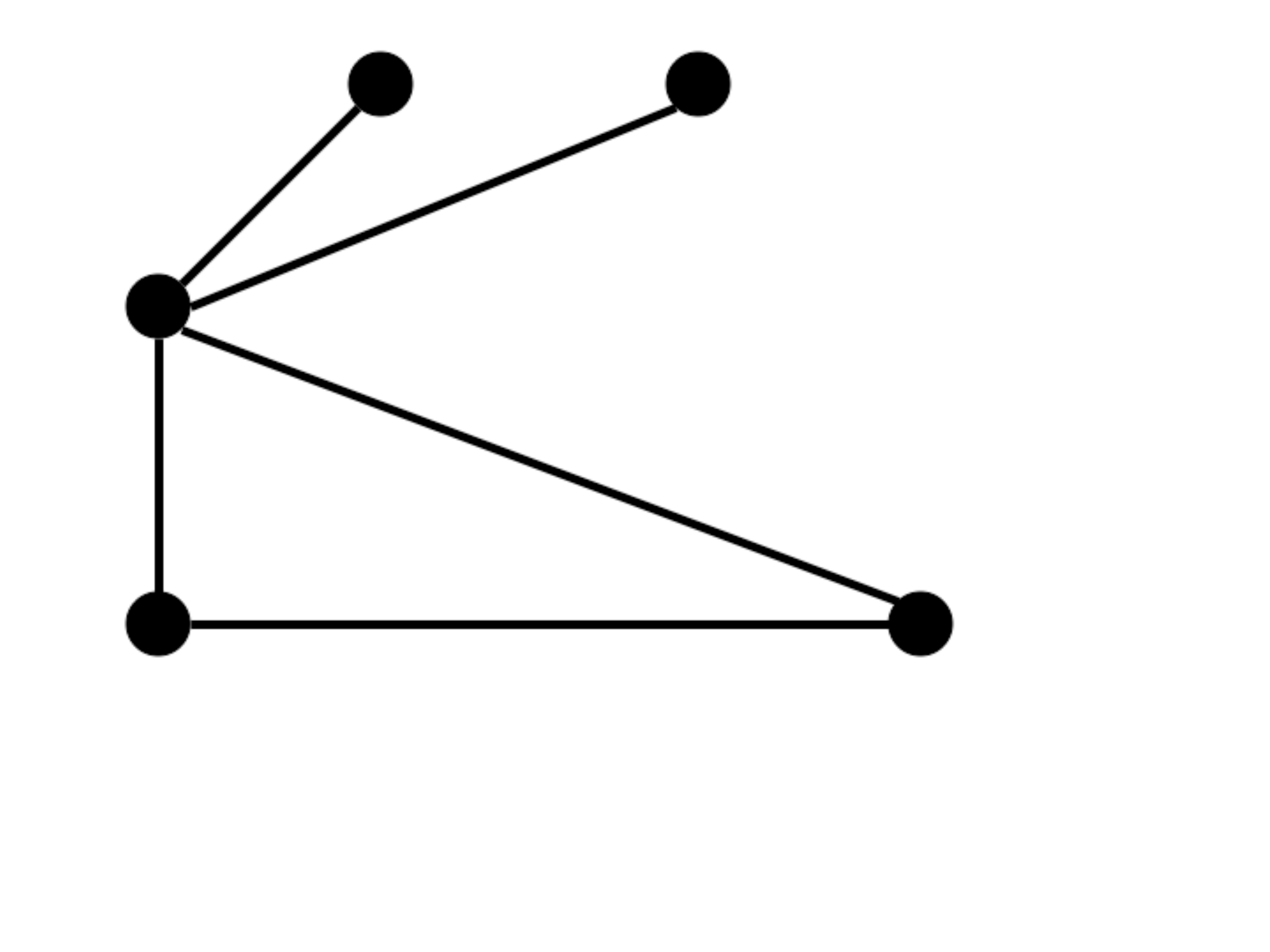} &
\includegraphics[width=0.6in]{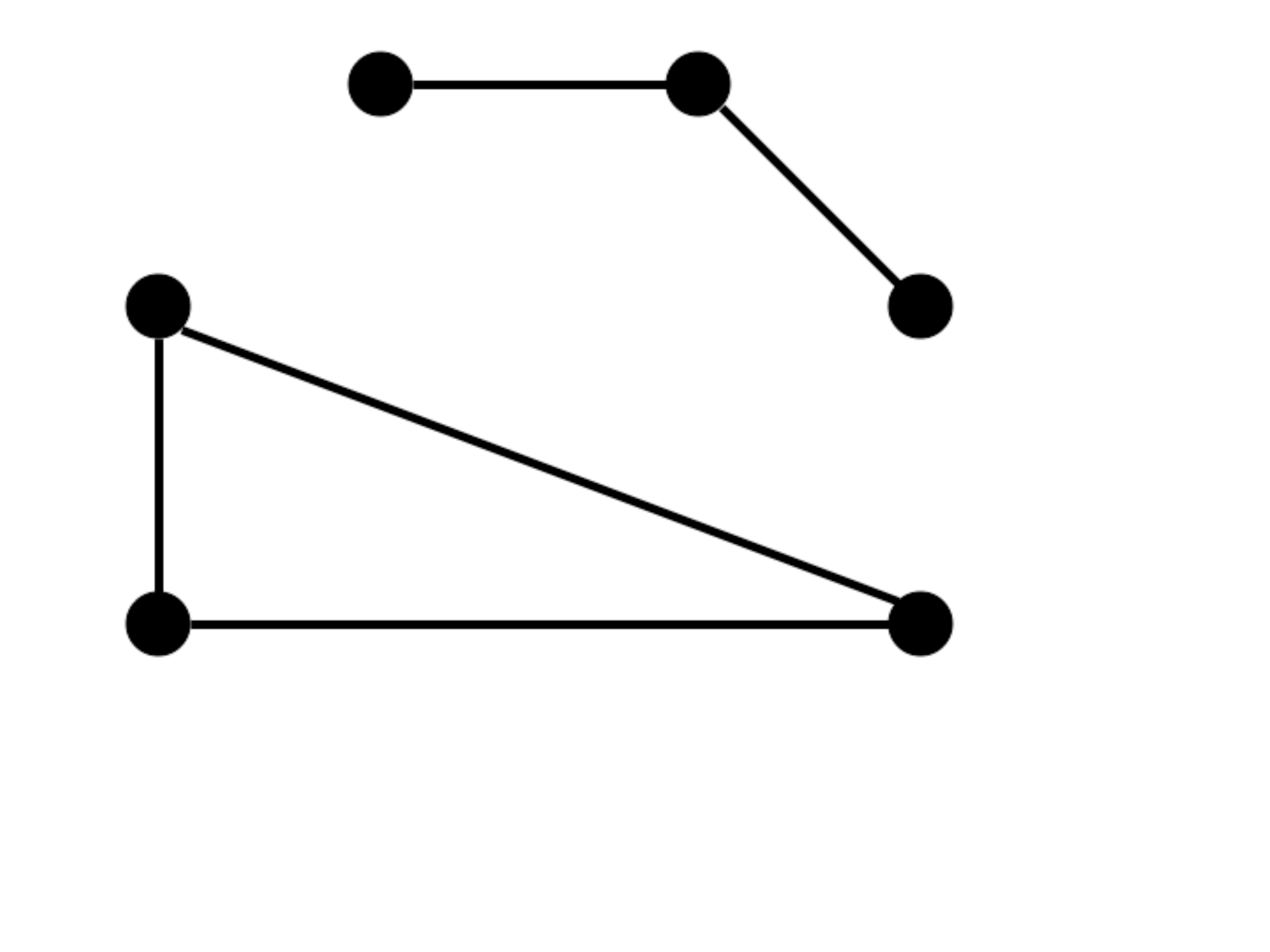} \\
\hline
\hline
$G_7^{(5)}$ & $G_8^{(5)}$ & $G_9^{(5)}$ & $G_{10}^{(5)}$ & $G_{11}^{(5)}$ & $G_{12}^{(5)}$ \\
\hline
\includegraphics[width=0.6in]{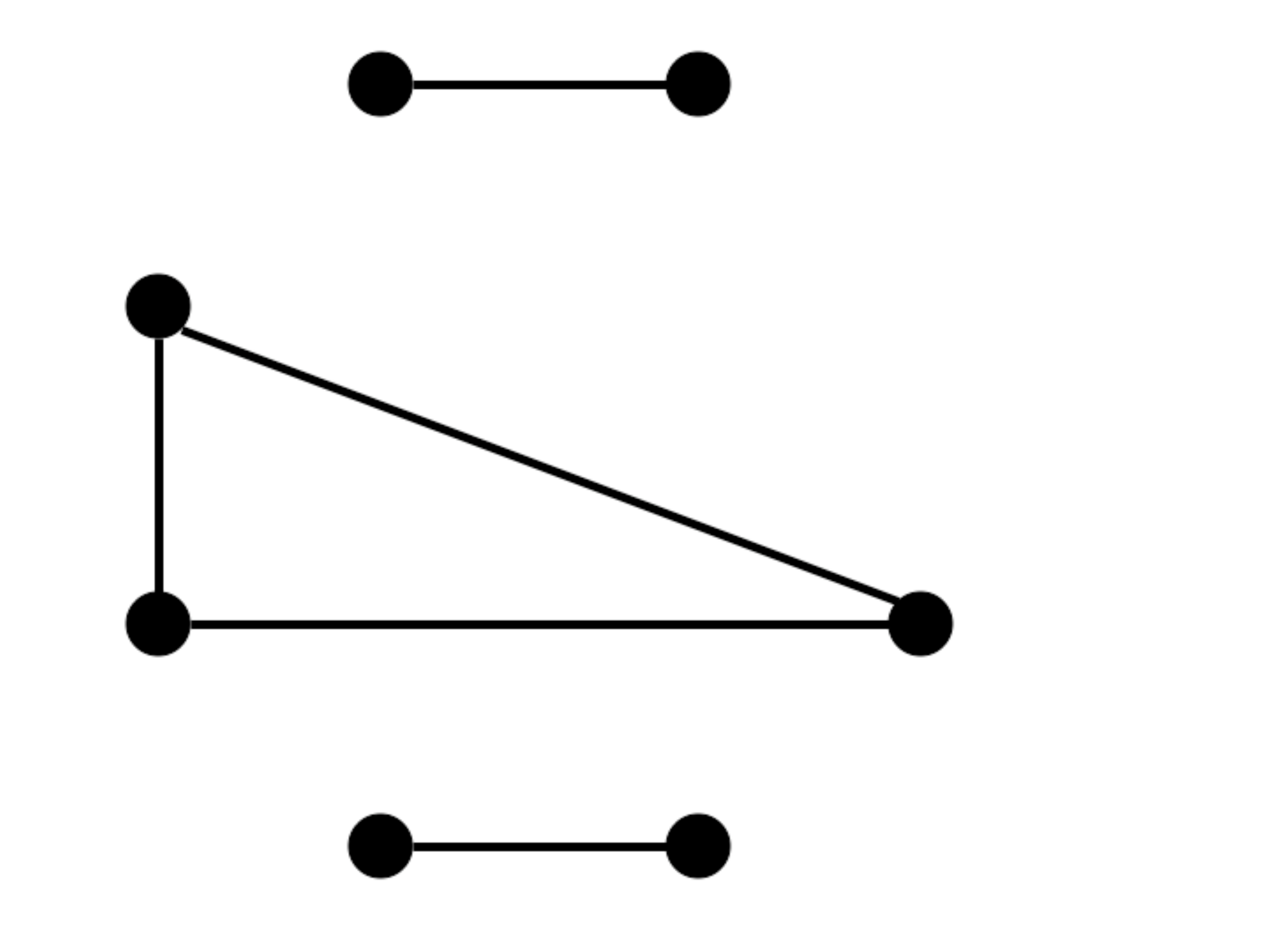} &
\includegraphics[width=0.6in]{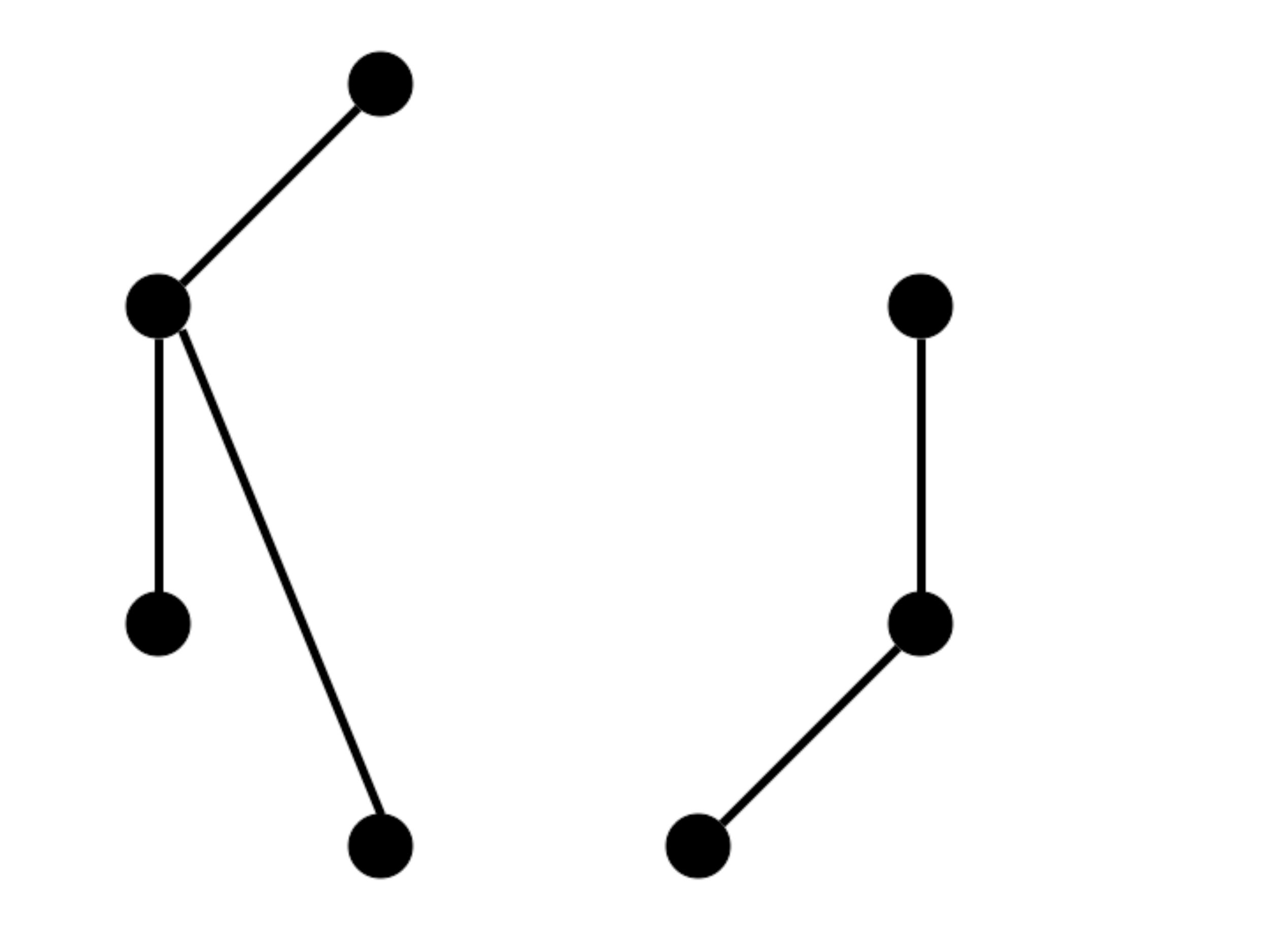} &
\includegraphics[width=0.6in]{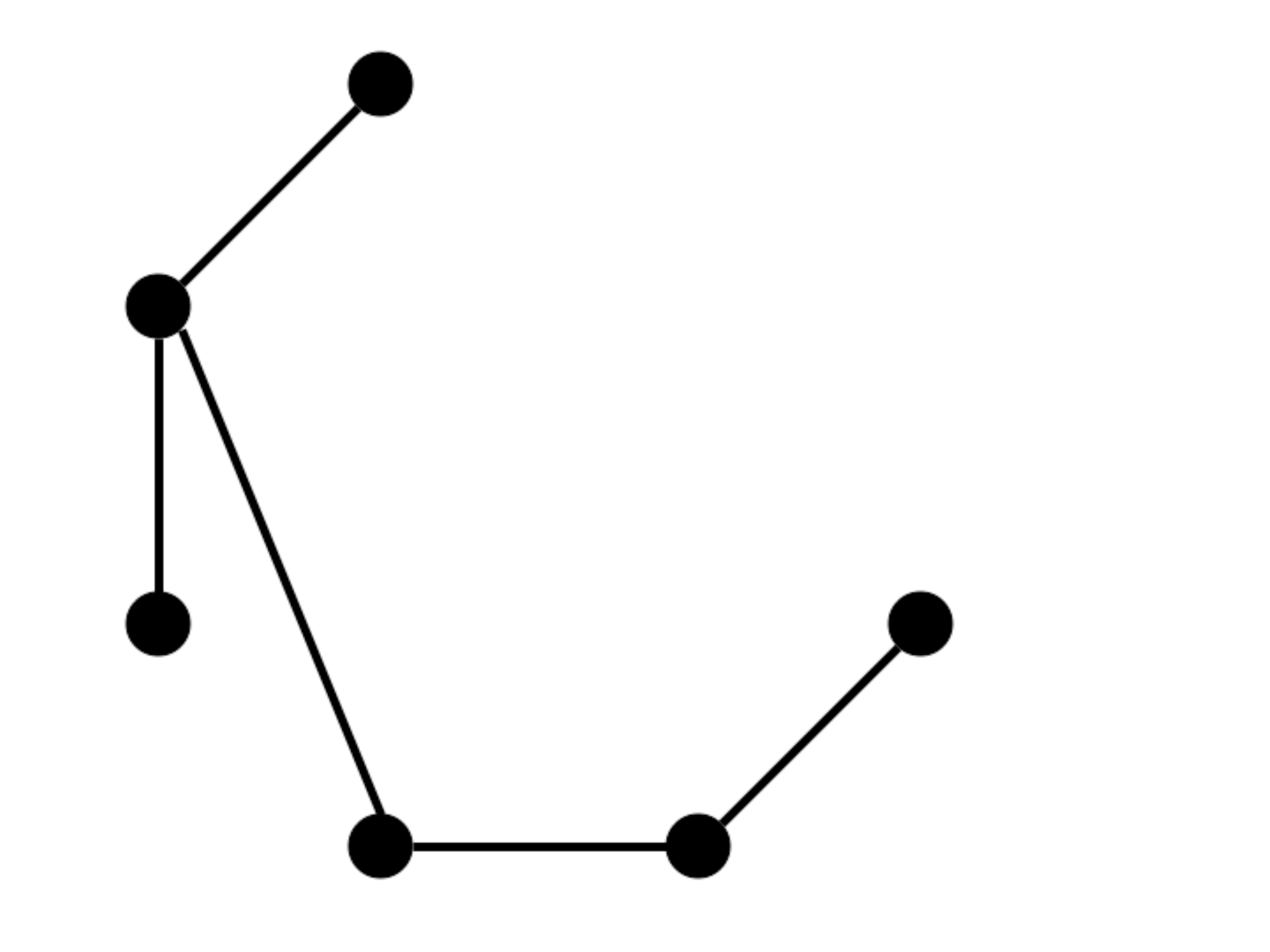} &
\includegraphics[width=0.6in]{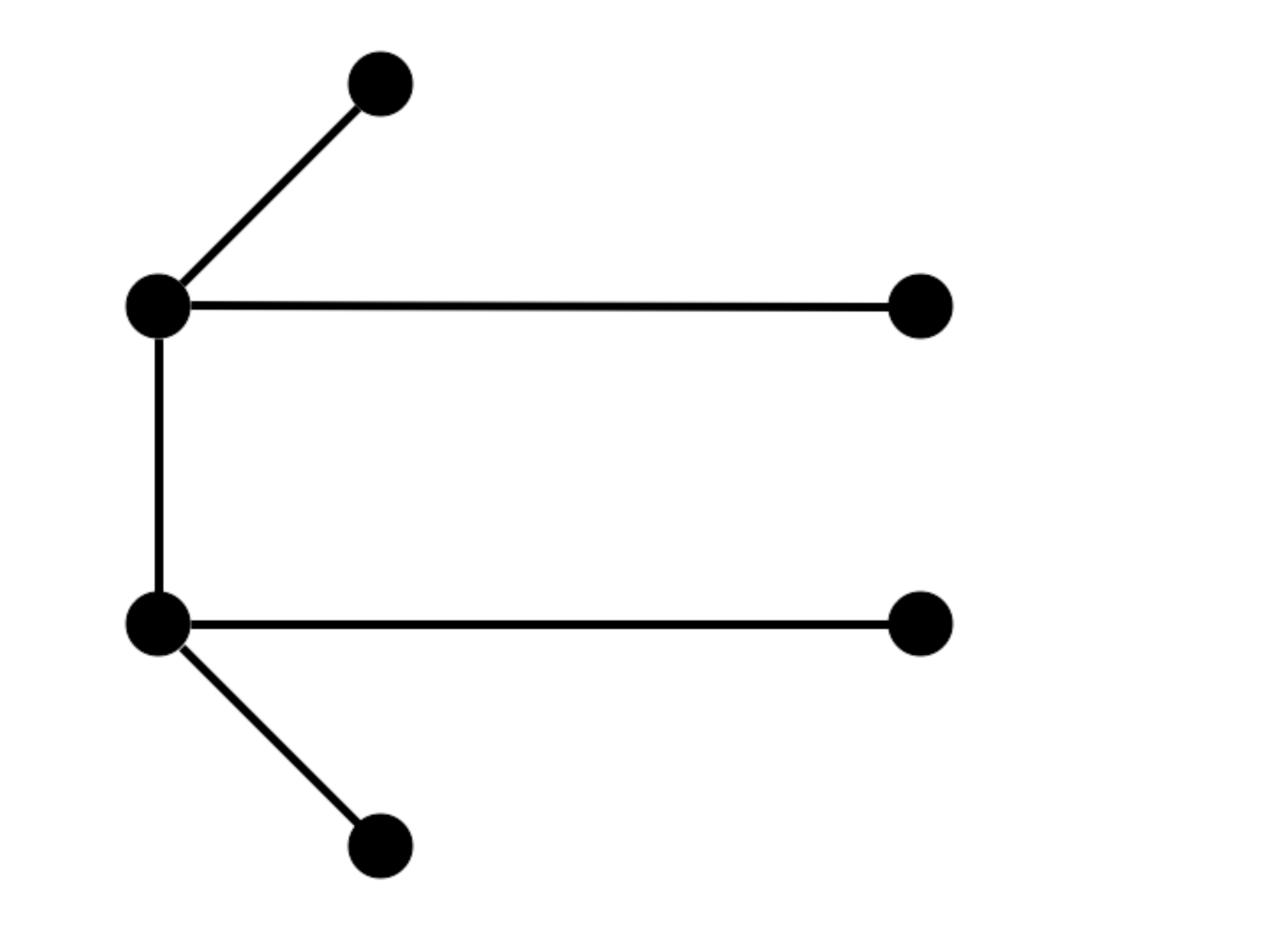} &
\includegraphics[width=0.6in]{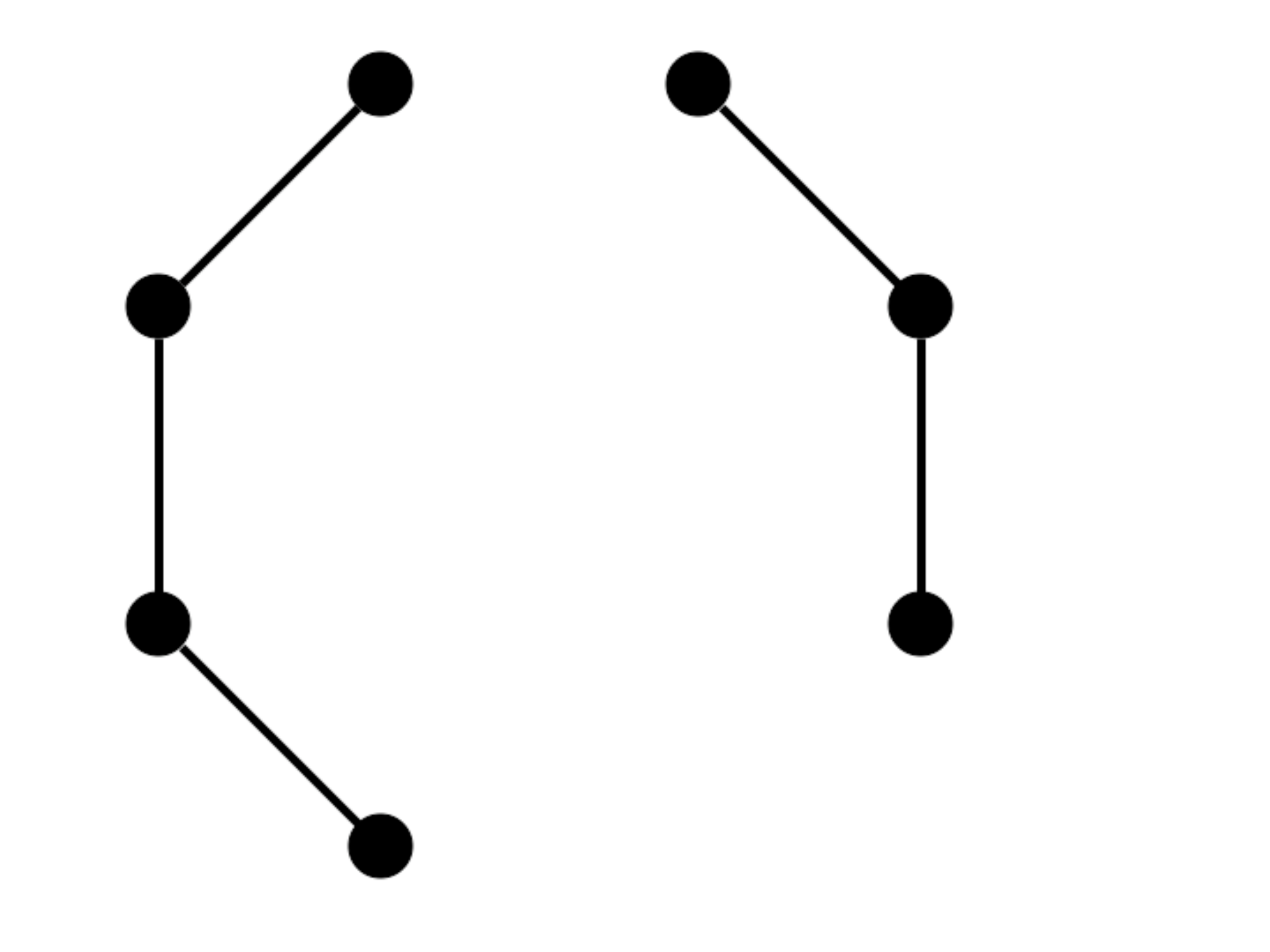} &
\includegraphics[width=0.6in]{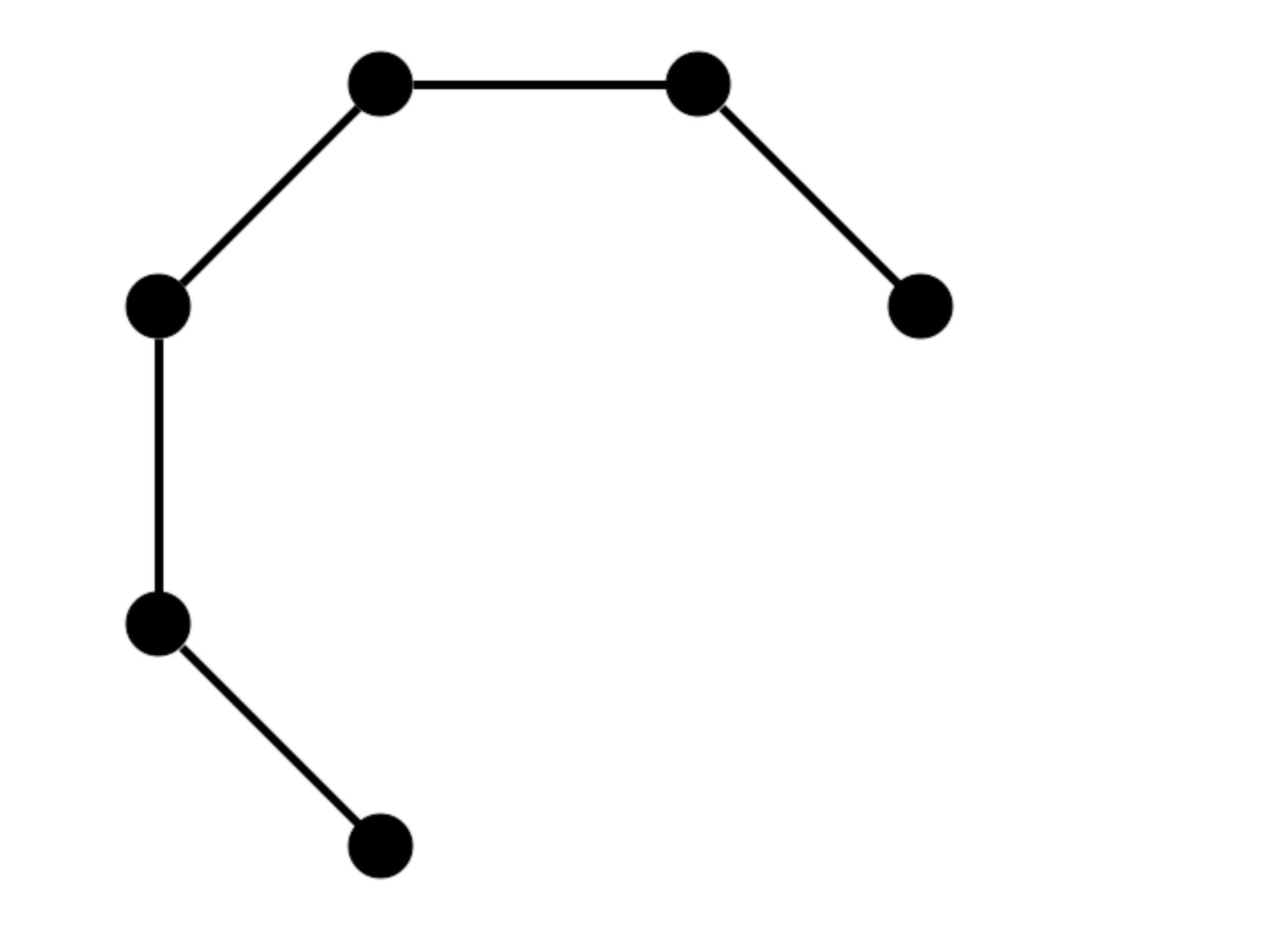} \\
\hline
\hline
$G_{13}^{(5)}$ & $G_{14}^{(5)}$ & $G_{15}^{(5)}$ & $G_{16}^{(5)}$ & $G_{17}^{(5)}$ & $G_{18}^{(5)}$ \\
\hline
\includegraphics[width=0.6in]{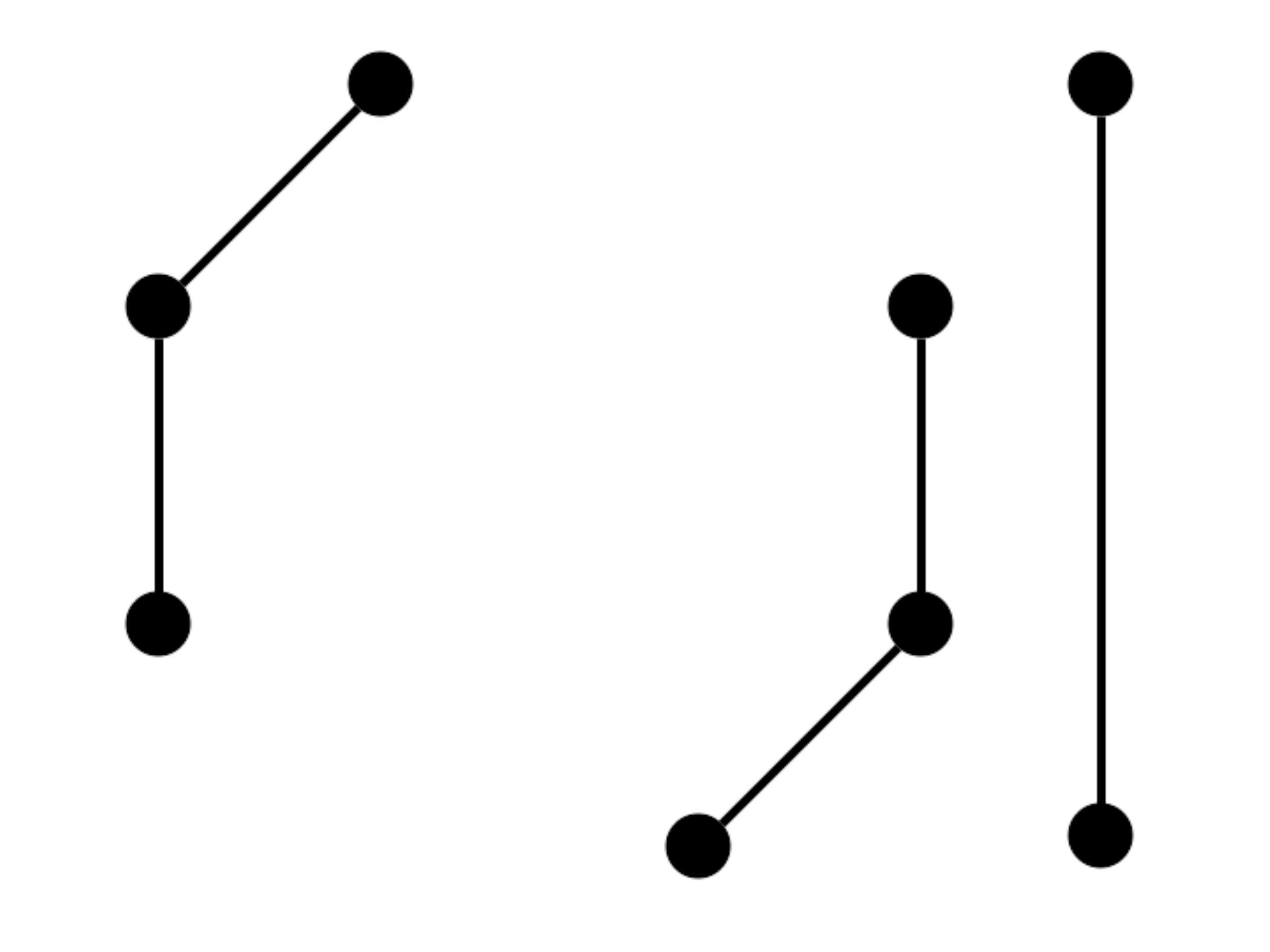} &
\includegraphics[width=0.6in]{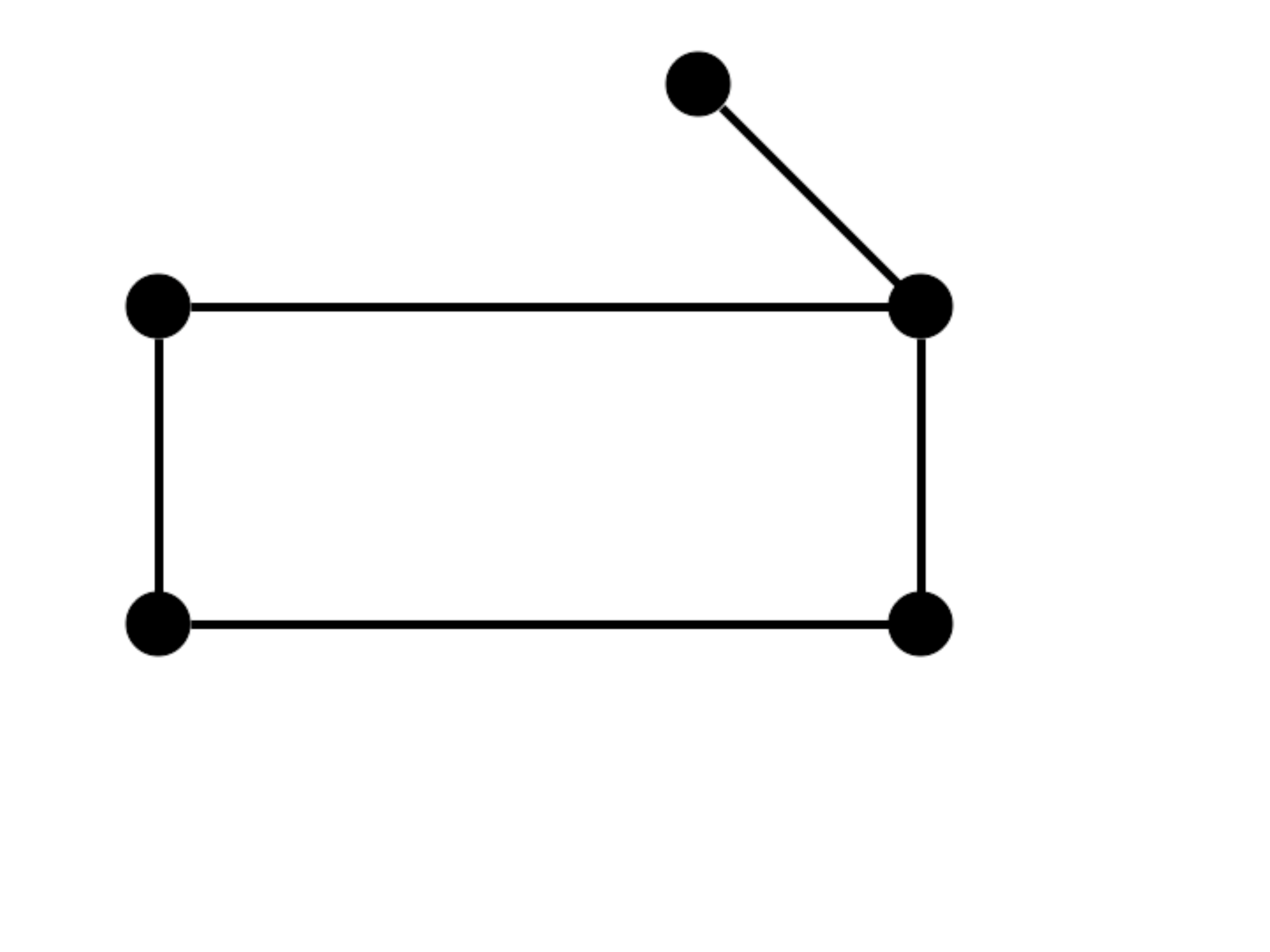} &
\includegraphics[width=0.6in]{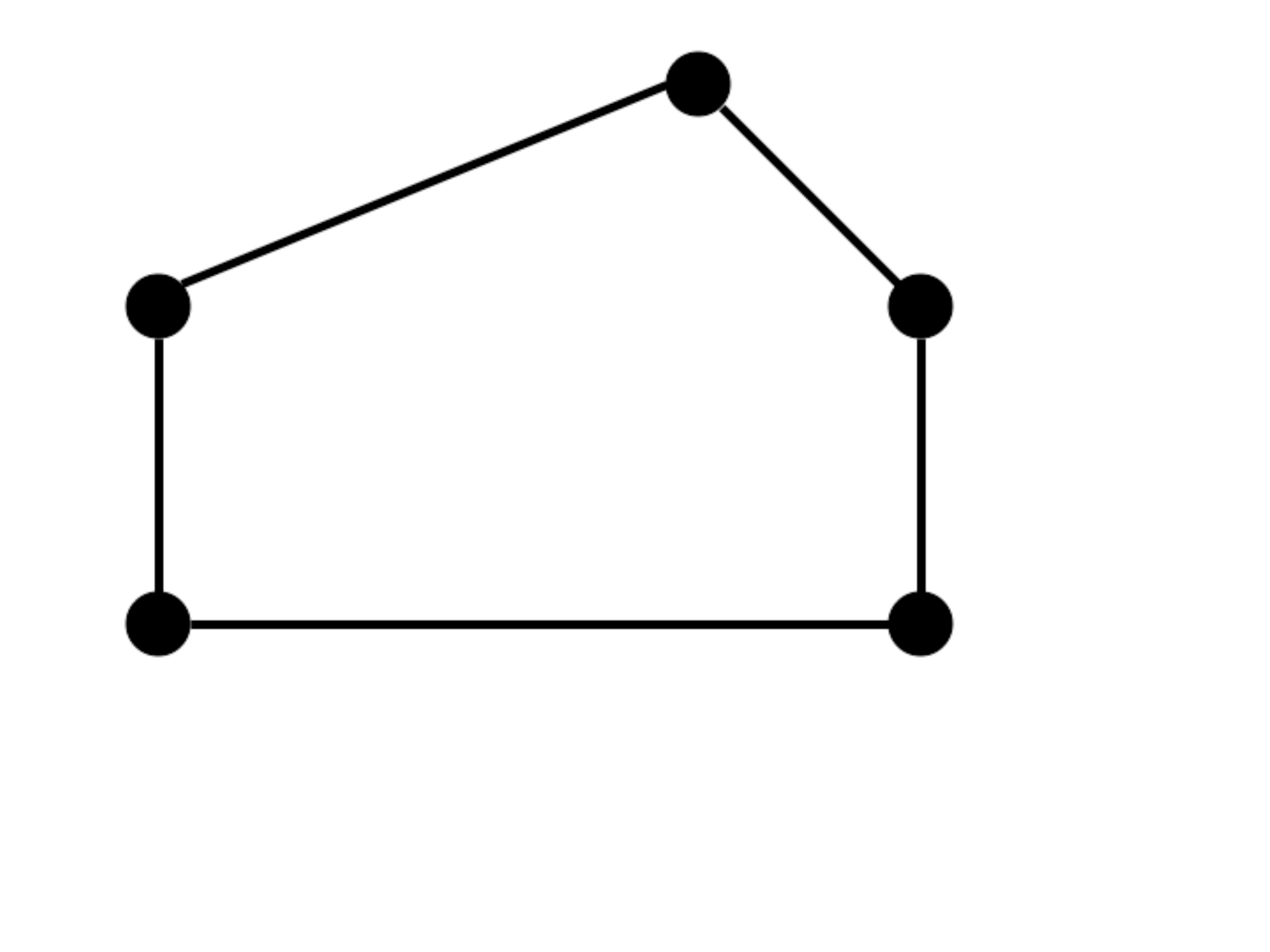} &
\includegraphics[width=0.6in]{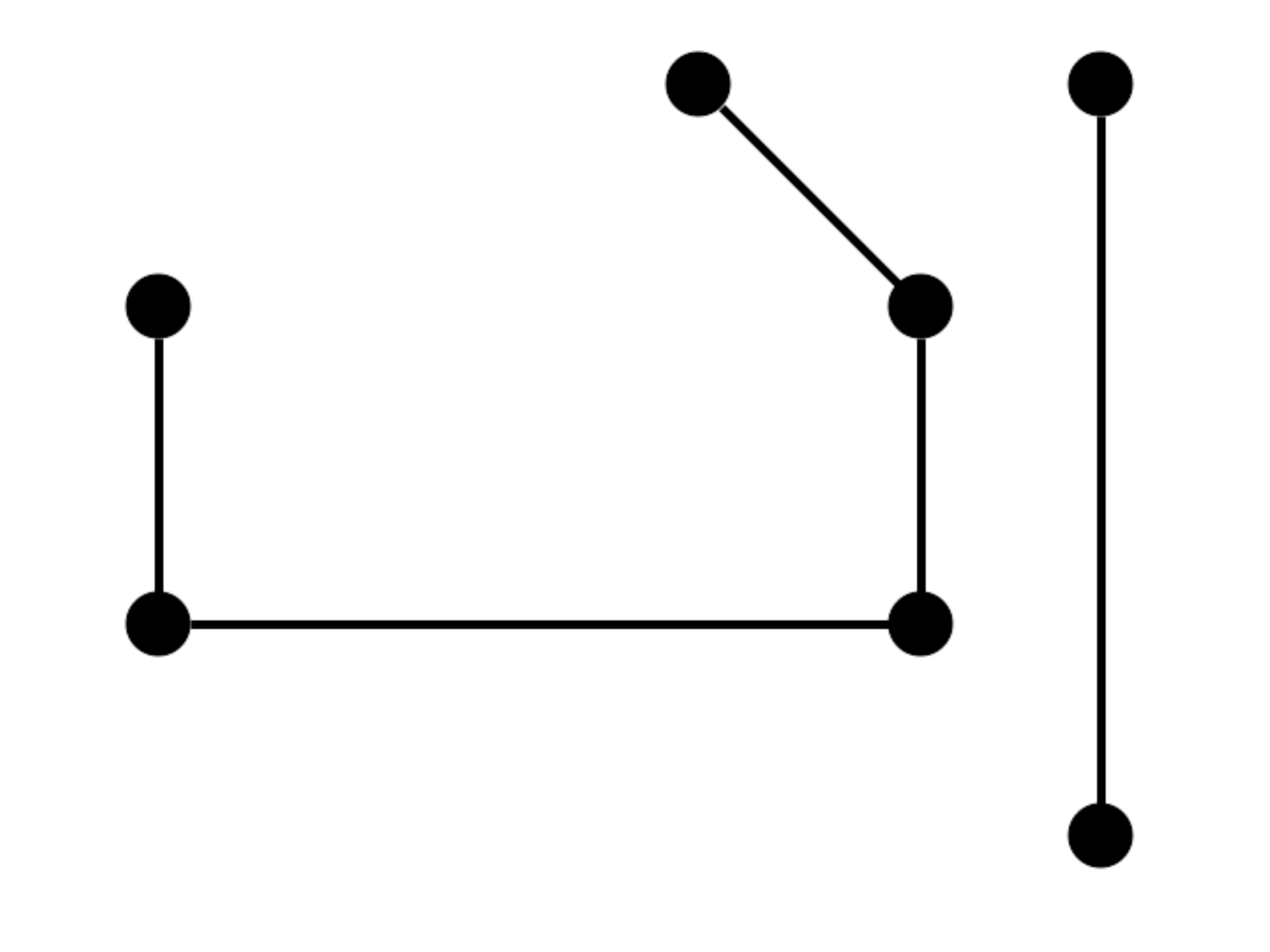} &
\includegraphics[width=0.6in]{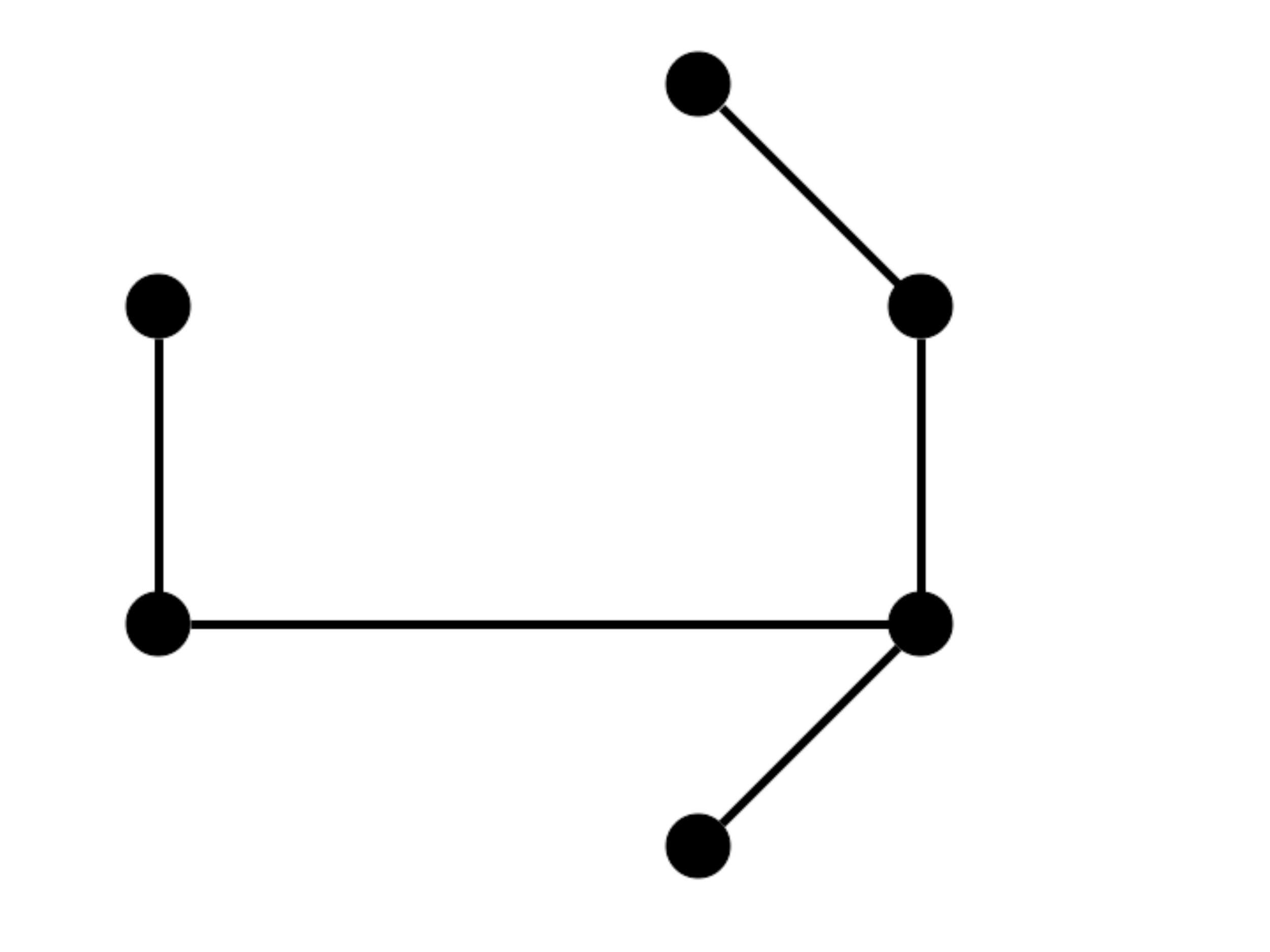} &
\includegraphics[width=0.6in]{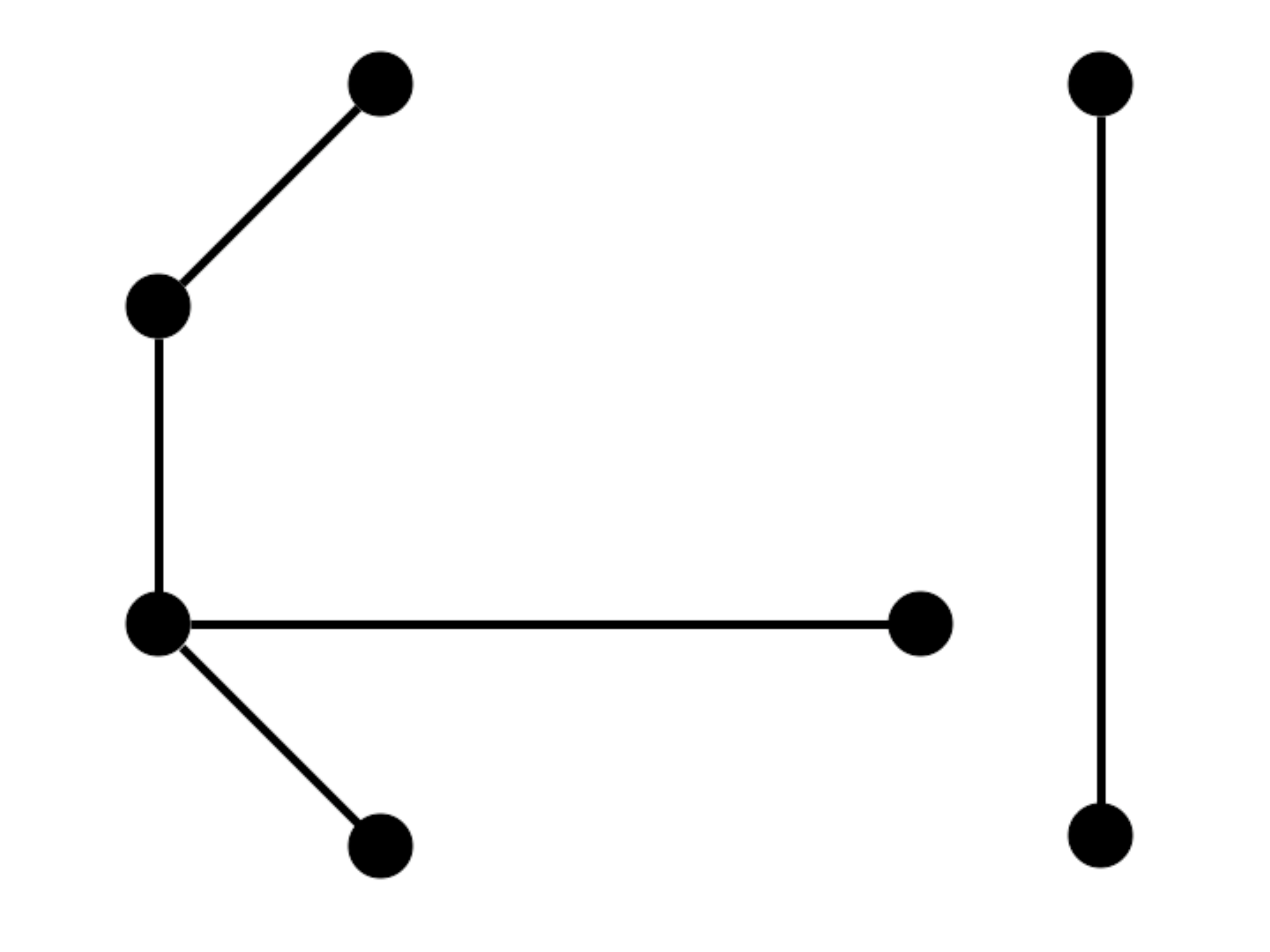} \\
\hline
\hline
$G_{19}^{(5)}$ & $G_{20}^{(5)}$ & $G_{21}^{(5)}$ & $G_{22}^{(5)}$ & $G_{23}^{(5)}$ & $G_{24}^{(5)}$ \\
\hline
\includegraphics[width=0.6in]{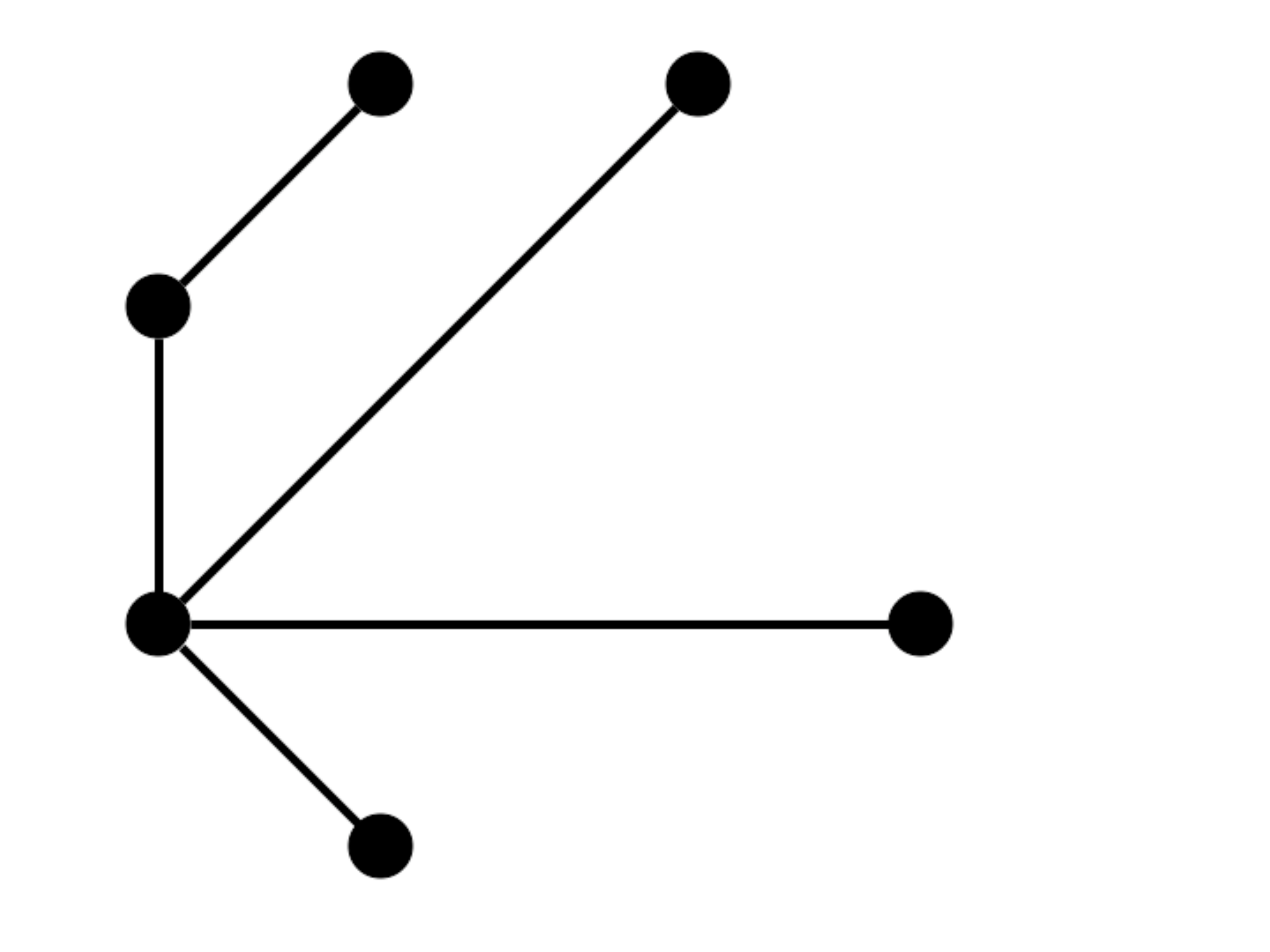} &
\includegraphics[width=0.6in]{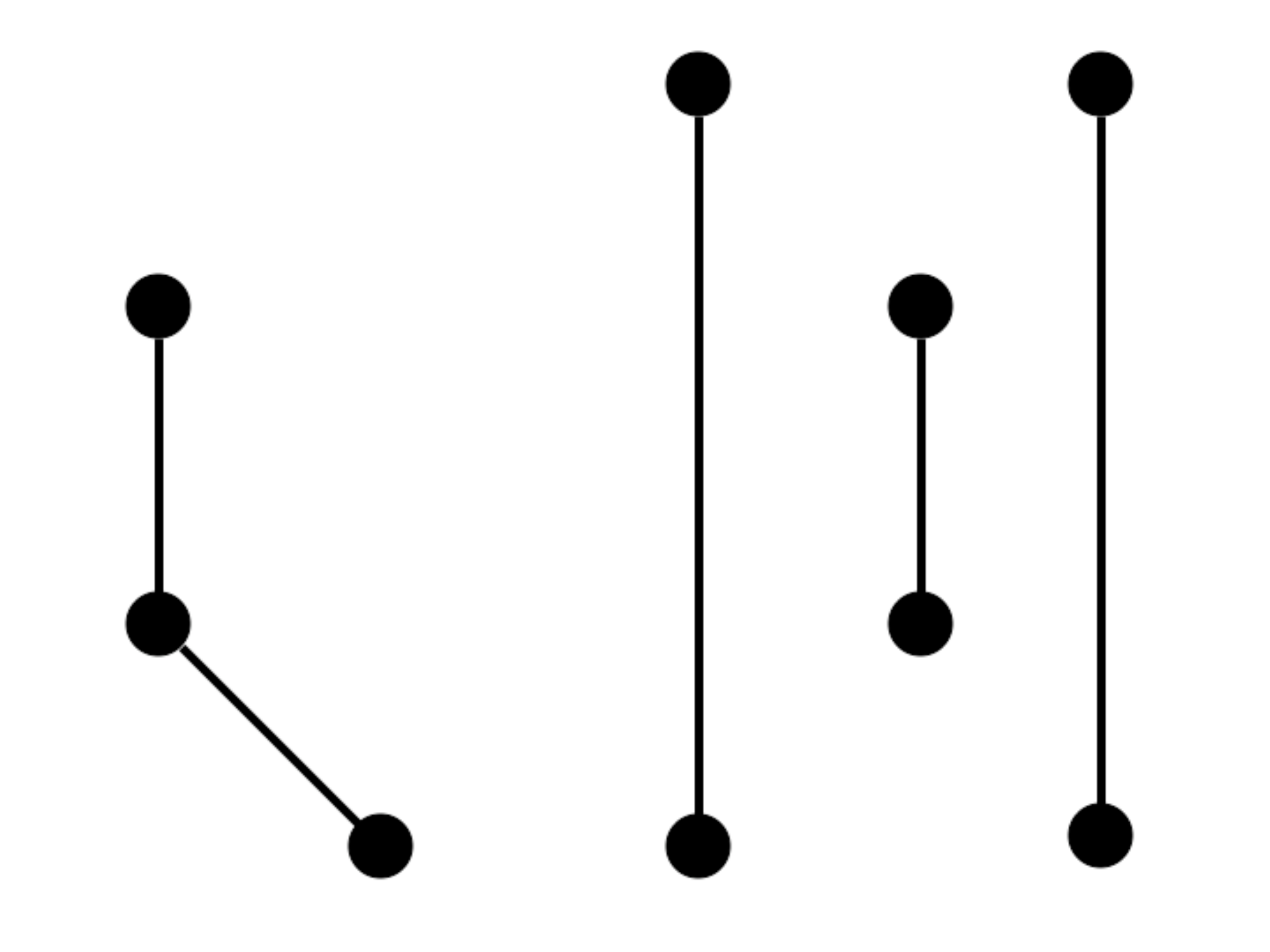} &
\includegraphics[width=0.6in]{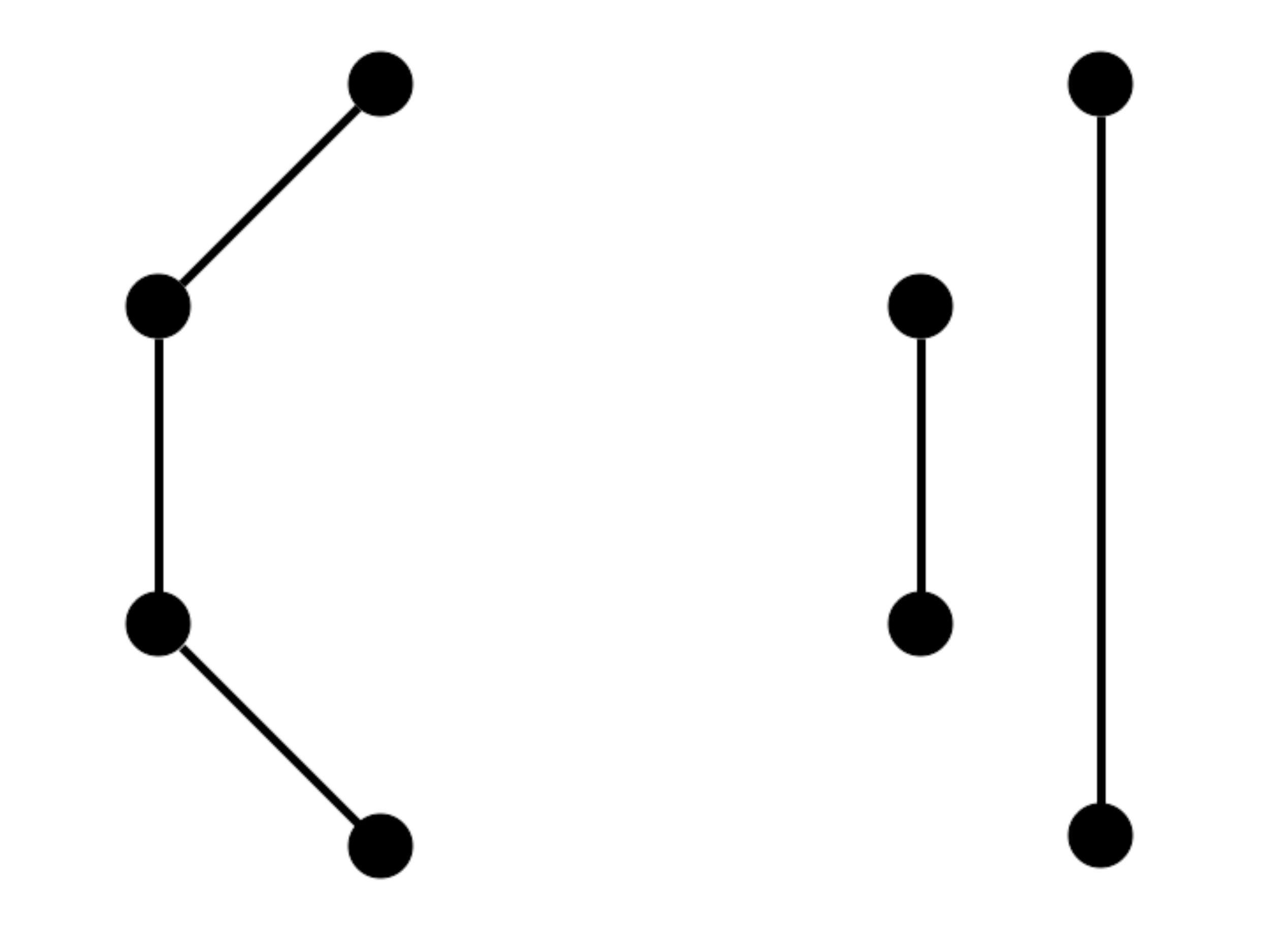} &
\includegraphics[width=0.6in]{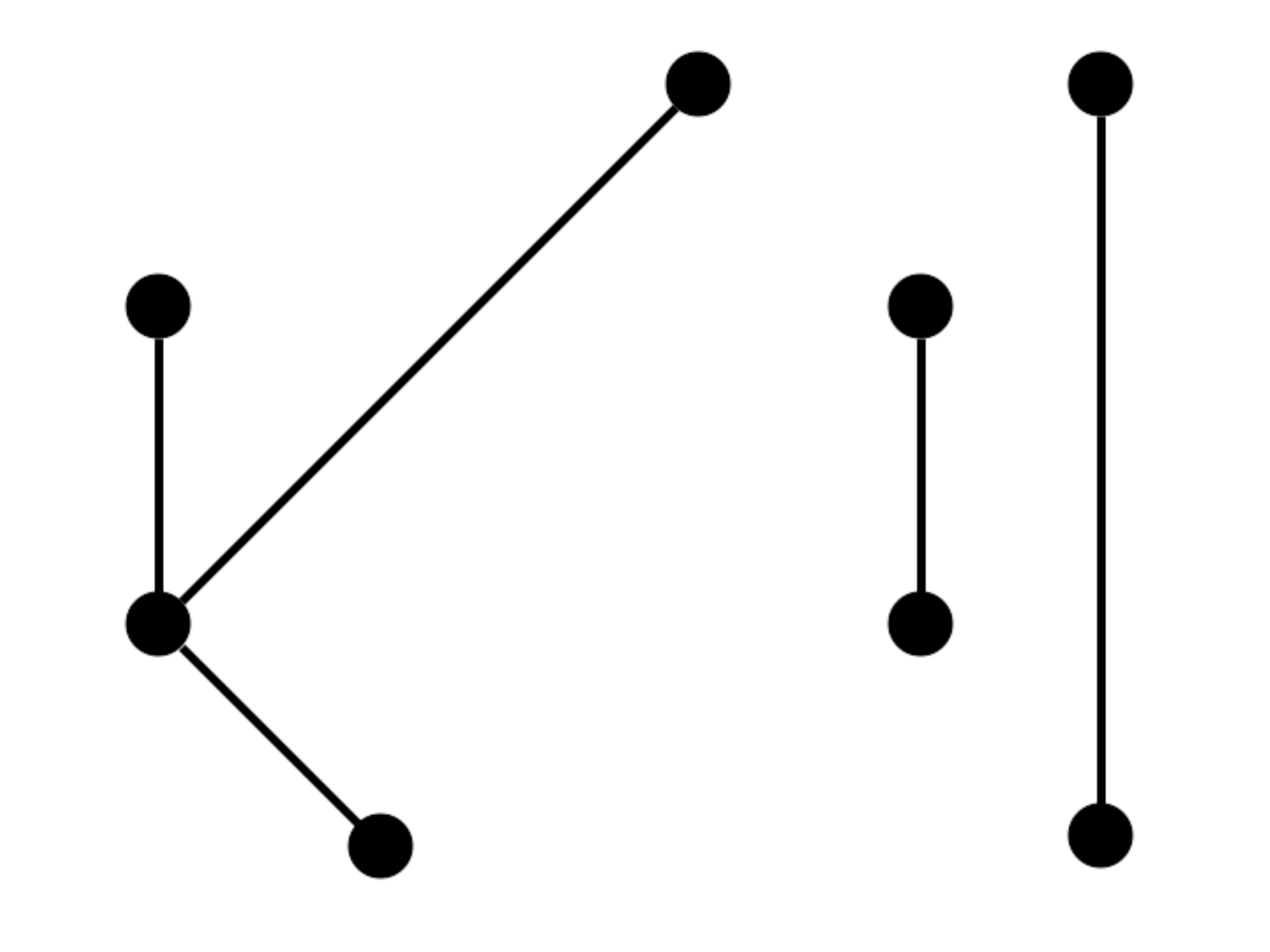} &
\includegraphics[width=0.6in]{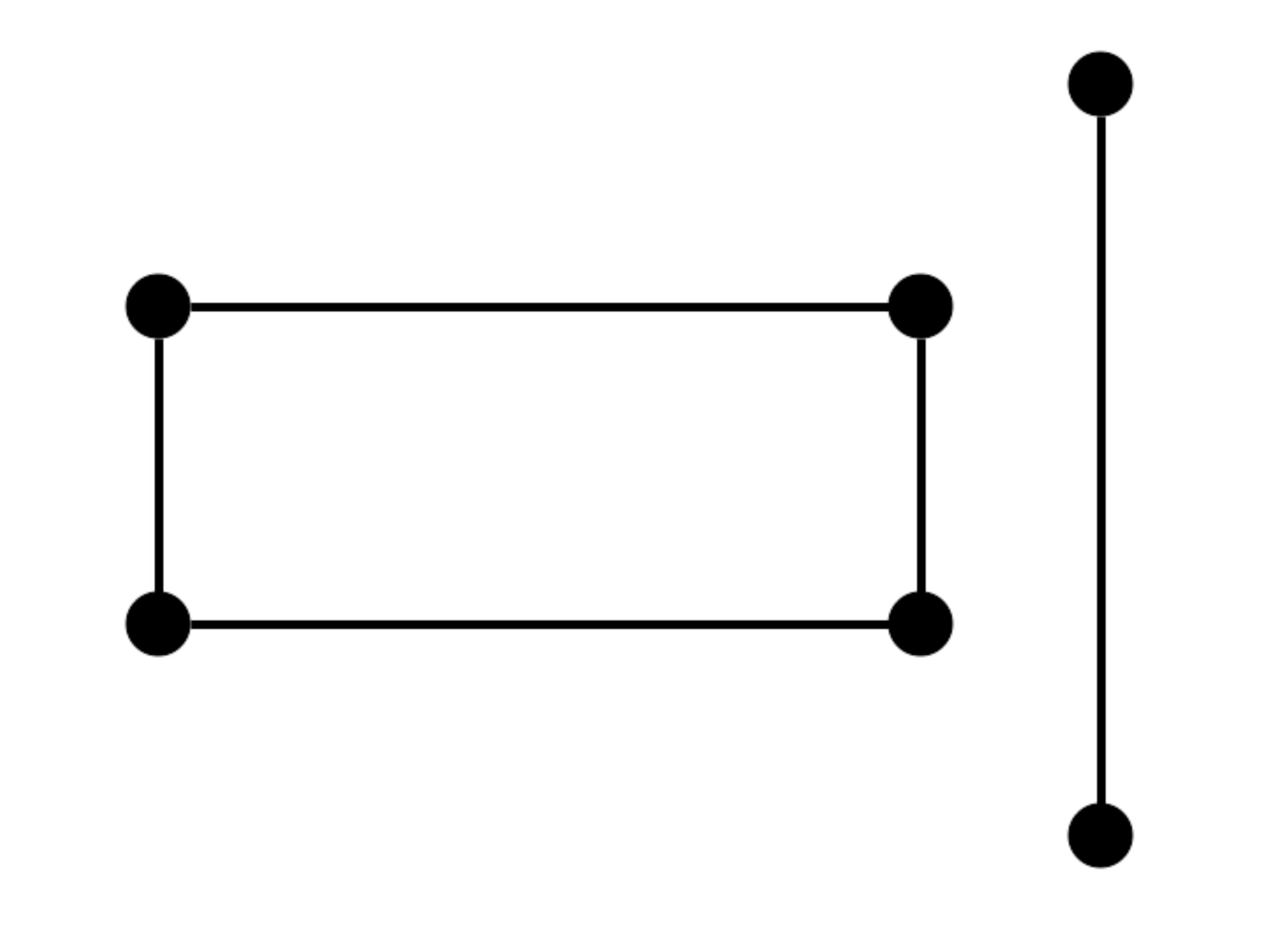} &
\includegraphics[width=0.6in]{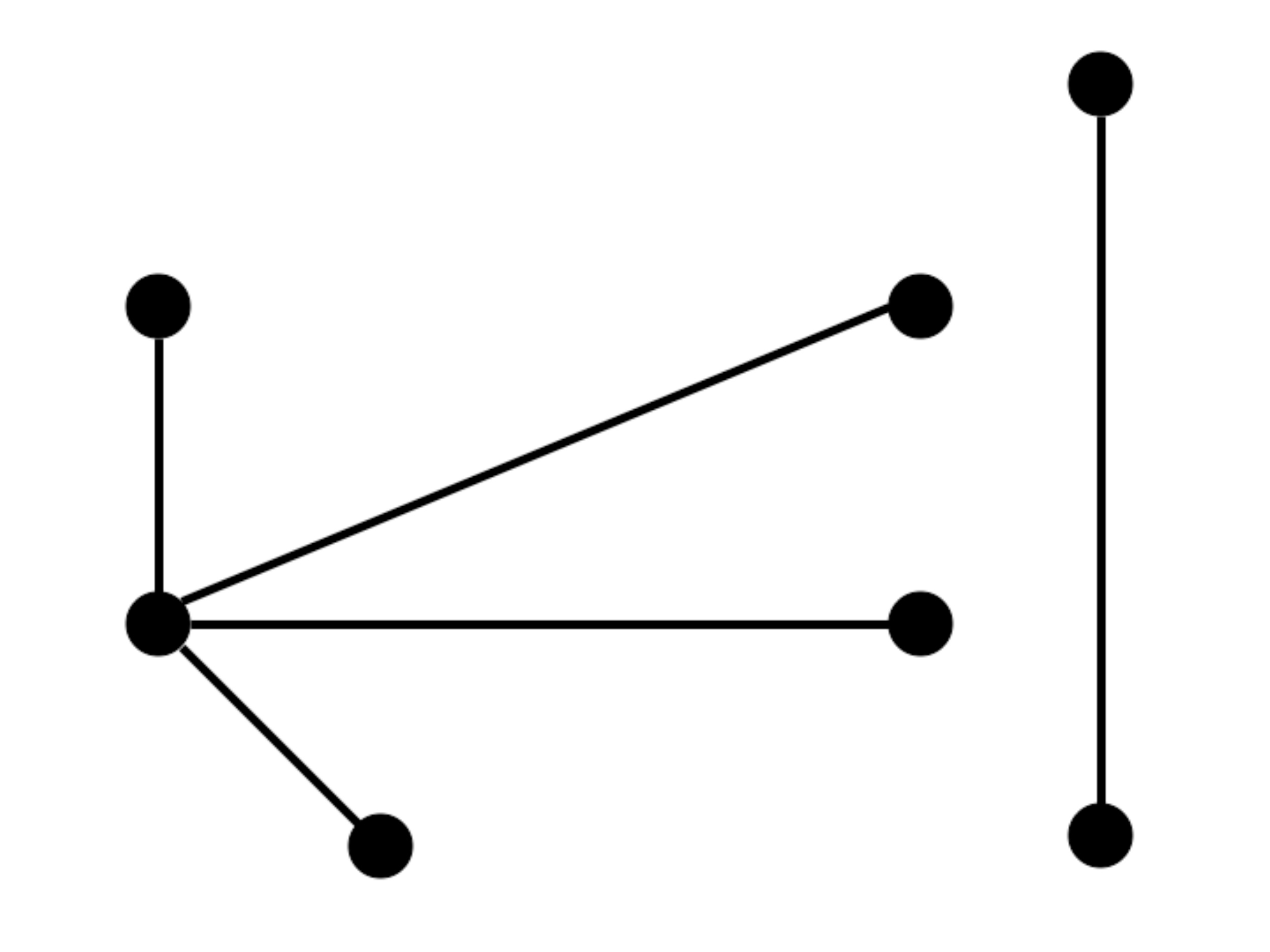} \\
\hline
\hline
$G_{25}^{(5)}$ & $G_{26}^{(5)}$ &  & &  &  \\
\hline
\includegraphics[width=0.6in]{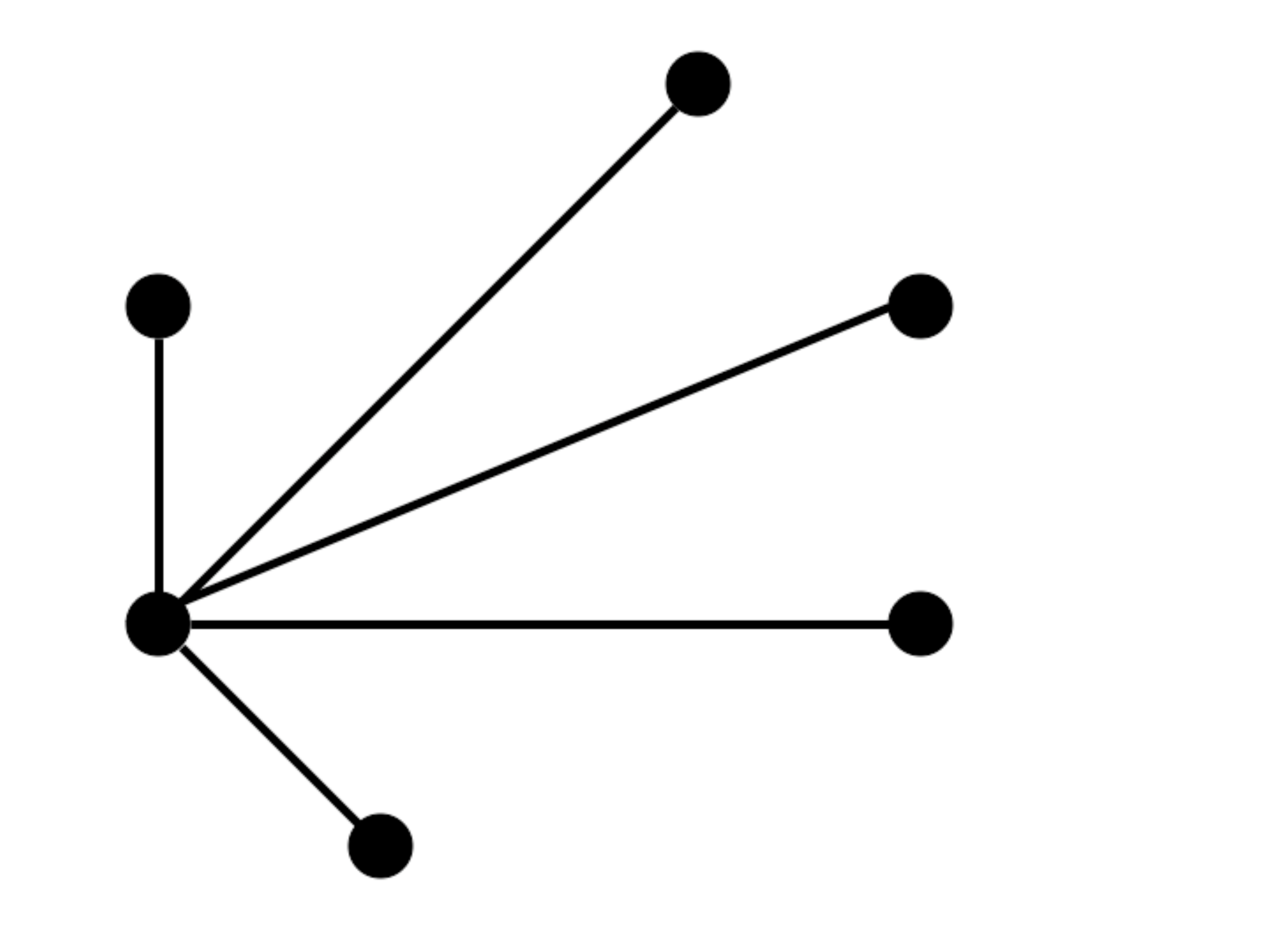} &
\includegraphics[width=0.6in]{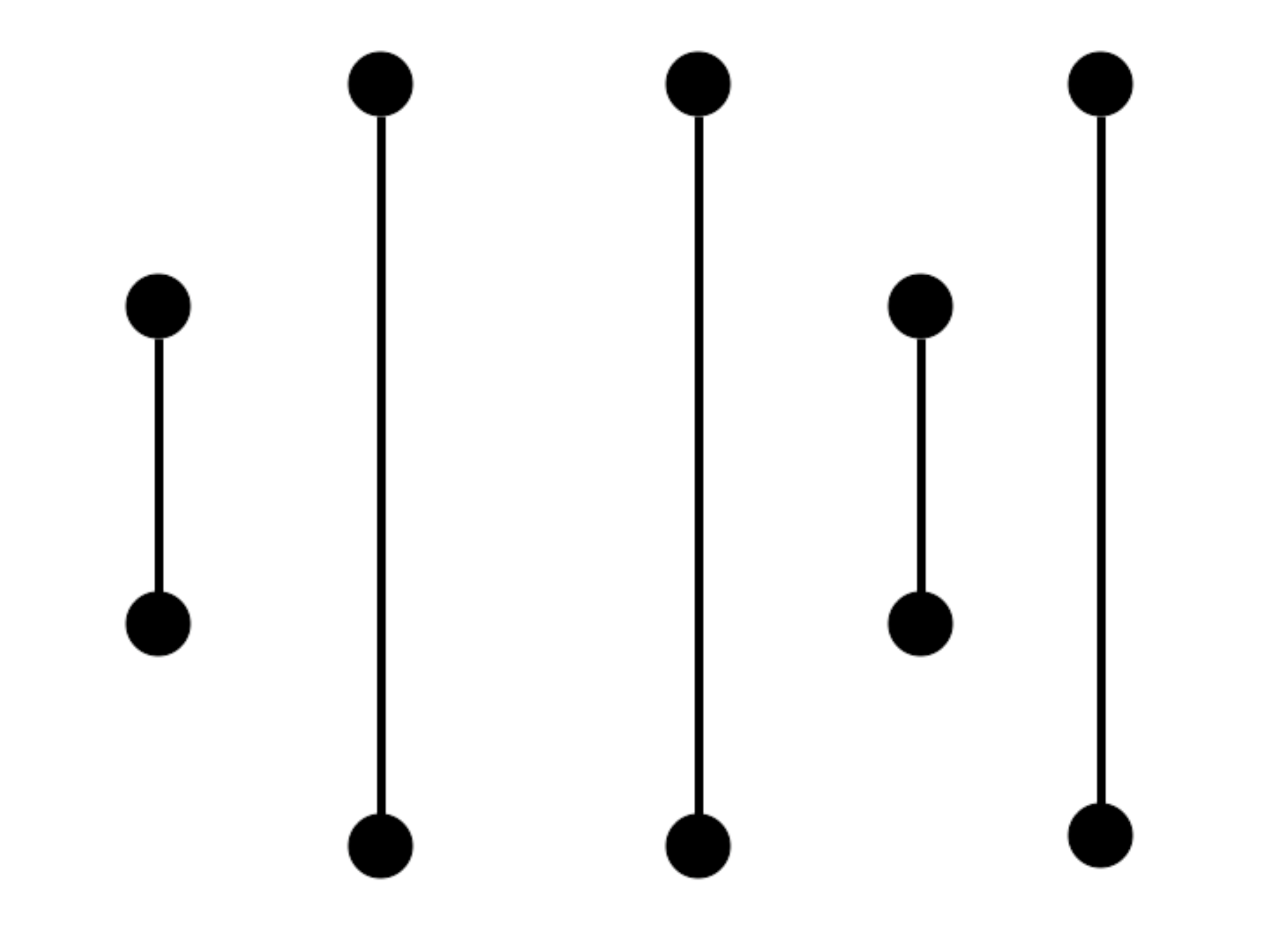} &
 & & & \\
\hline
\end{longtable}
}

\end{article}
\end{document}